\documentclass[12pt,leqno]{article}
\usepackage{times}
\usepackage{a4wide}
\usepackage{amssymb}
\usepackage{amsmath}
\usepackage{theorem}
\usepackage{enumerate}
\usepackage{amsmath}
\usepackage{amscd}
\usepackage{amssymb}

\theorembodyfont{\sl}
\newtheorem{subthm}[subsection]{Theorem}
\newtheorem{subprop}[subsection]{Proposition}

\theorembodyfont{\rmfamily}
\newtheorem{subrem}[subsection]{Remark}

\numberwithin{equation}{subsection}

\makeatletter \newcommand{\subsectionrunin} {\@startsection
{subsection}{1}{\z@ }{-3.5ex \@plus -1ex \@minus -.2ex}
{-1.5ex}{\normalfont \normalsize \bfseries }} \makeatother

\newenvironment{prf}[1]{\trivlist
\item[\hskip \labelsep{\it
#1.\hspace*{.3em}}]}{~\hspace{\fill}~$\square$\endtrivlist}
\newenvironment{proof}{\begin{prf}{\bf Proof}}{\end{prf}}

\newcommand{\CC}{\mathbb{C}}
\newcommand{\QQ}{\mathbb{Q}}
\newcommand{\FF}{\mathbb{F}}
\newcommand{\ZZ}{\mathbb{Z}}
\newcommand{\PP}{\mathbb{P}}
\newcommand{\EE}{\mathbb{E}}
\newcommand{\GG}{\mathbb{G}}
\newcommand{\TT}{\mathbb{T}}
\newcommand{\HH}{\mathbb{H}}

\newcommand{\Qbar}{{\overline{\QQ}}}
\newcommand{\Fbar}{{\overline{\FF}}}
\DeclareMathOperator{\Spec}{Spec}
\DeclareMathOperator{\Cot}{Cot}
\newcommand{\Katz}{\mathrm{Katz}}

\newcommand{\rH}{{\mathrm{H}}}
\newcommand{\rR}{{\mathrm{R}}}

\DeclareMathOperator{\dash_exp}{-exp}

\newcommand{\calO}{\mathcal{O}}

\newcommand{\ol}{\overline}
\newcommand{\ul}{\underline}
\newcommand{\lto}{\longrightarrow}
\newcommand{\sm}{\mathrm{sm}}
\newcommand{\eps}{\varepsilon}

\newcommand{\End}{{\rm End}}
\newcommand{\Aut}{{\rm Aut}}
\newcommand{\Hom}{{\rm Hom}}
\newcommand{\GL}{\mathrm{GL}}
\newcommand{\SL}{\mathrm{SL}}
\newcommand{\Gm}[2]{{\GG_{\mathrm{m}{{#1}}}^{{#2}}}}
\newcommand{\ld}{{\langle}}
\newcommand{\rd}{{\rangle}}
\newcommand{\cusps}{\mathrm{cusps}}
\DeclareMathOperator{\Sym}{Sym}
\newcommand{\parabolic}{\mathrm{par}}
\newcommand{\compact}{\mathrm{c}} 
\DeclareMathOperator{\Image}{Image}
\newcommand{\et}{\mathrm{et}}
\DeclareMathOperator{\Gal}{Gal}
\DeclareMathOperator{\Fil}{Fil}
\DeclareMathOperator{\Gr}{Gr}
\newcommand{\deRham}{\mathrm{dR}} 

\DeclareMathOperator{\Rep}{Rep}
\DeclareMathOperator{\MF}{MF}
\DeclareMathOperator{\modulo}{mod}
\DeclareMathOperator{\Frob}{Frob}

\newtheorem{subsubnot}[subsubsection]{Notation}
\newtheorem{subsubsetting}[subsubsection]{Setting}
\newtheorem{subsubcor}[subsubsection]{Corollary}
\newtheorem{subsubprop}[subsubsection]{Proposition}

\newcommand{\cS}{\mathcal{S}}
\newcommand{\cR}{\mathcal{R}}
\newcommand{\cO}{\mathcal{O}}
\newcommand{\cA}{\mathcal{A}}
\newcommand{\cH}{\mathcal{H}}
\newcommand{\fm}{\mathfrak{m}}
\newcommand{\fP}{\mathfrak{P}}

\newcommand{\ttt}[1]{{\tt #1}}

\newcommand{\id}{\mathrm{id}}
\newcommand{\Bound}{{\rm Bound}}
\newcommand{\Ind}{{\rm Ind}}
\newcommand{\PGL}{\mathrm{PGL}}

\newcommand{\ksf}[4]{S_{#1} (\Gamma_{#2} (#3), {#4})_\Katz'}
\newcommand{\ksfc}[5]{S_{#1} (\Gamma_{#2} (#3), {#4}, {#5})_\Katz' }
\newcommand{\Sf}[4]{S_{#1} (\Gamma_{#2} (#3), {#4})}
\newcommand{\Sfc}[5]{S_{#1} (\Gamma_{#2} (#3), {#4}, {#5}) }
\newcommand{\ms}[3]{\cS_{#1} (\Gamma_{#2} (#3))}
\newcommand{\msc}[4]{\cS_{#1} (\Gamma_{#2} (#3), {#4}) }

\newcommand{\comp}[1]{\begin{small}\begin{quote}{\tt #1}\end{quote}\end{small}}

\begin{document}

\title{Comparison of integral structures on spaces of modular forms of
weight two, and computation of spaces of forms mod 2 of weight one\\
\Large{with appendices by Jean-Fran{\c c}ois Mestre and Gabor Wiese}}

\author{Bas Edixhoven\footnote{Partially supported by the European
Research Training Network Contract HPRN-CT-2000-00120 ``arithmetic
algebraic geometry''. Address: Mathematisch Instituut, Universiteit
Leiden, Postbus 9512, 2300\ RA\ \ Leiden, The Netherlands,
edix@math.leidenuniv.nl}}

\maketitle

\setcounter{tocdepth}{1}
\tableofcontents

\section{Introduction} 
For $N\geq1$ and $k$ integers, let $S_k(\Gamma_0(N),\CC)$ denote the
$\CC$-vector space of holomorphic cusp forms on the congruence
subgroup $\Gamma_0(N)$ of~$\SL_2(\ZZ)$ (for a definition, see for
example the book containing~\cite{Edix5}).  To an element $f$ of
some $S_k(\Gamma_0(N),\CC)$ we can associate its $q$-expansion
$\sum_{n\geq1}a_n(f)q^n$ in~$\CC[[q]]$, by noting that~$f$, as a
function on the complex upper half plane, is invariant under
translation by all integers, and hence is a power series in $q\colon
z\mapsto\exp(2\pi iz)$. So we have $q$-expansion maps:
$$
q\dash_exp\colon S_k(\Gamma_0(N),\CC) \lto \CC[[q]], \qquad
f\mapsto \sum_{n\geq1}a_n(f)q^n, 
$$
which are injective. For $\Gamma$ any subgroup of finite index of
$\SL_2(\ZZ)$ and $k$ an integer we have an analogous space
$S_k(\Gamma,\CC)$. 

The $\CC$-vector spaces $S_k(\Gamma_0(N),\CC)$ have natural
$\QQ$-structures, i.e., sub $\QQ$-vector spaces $S_k(\Gamma_0(N),\QQ)$
such that the natural maps from $\CC\otimes S_k(\Gamma_0(N),\QQ)$ to
$S_k(\Gamma_0(N),\CC)$ are isomorphisms. One way to define these
$\QQ$-structures is to take the sub $\QQ$-vector space of $f$ in
$S_k(\Gamma_0(N),\CC)$ whose $q$-expansion is in~$\QQ[[q]]$. Note that
it is not a priori clear that these subspaces indeed are
$\QQ$-structures.  A second way to define the $\QQ$-structures is to
let $S_k(\Gamma_0(N),\QQ)$ be the space of global sections of the
invertible sheaf of modules~$\ul{\omega}^{\otimes k}(-\cusps)$ on the
$\QQ$-stack of generalized elliptic curves over $\QQ$-schemes with a
$\Gamma_0(N)$-level structure. The two $\QQ$-structures agree, but
this is not a complete triviality; see below for more details.

Loosely speaking, the two $\QQ$-structures come from $q$-expansions
and algebraic geometry over~$\QQ$. These two methods also give
$\ZZ$-structures. The two $\ZZ$-structures do not always coincide.

The first aim of this article is not to study in general the
differences between the two kinds of $\ZZ$-structures, although this
is an interesting problem. We restrict ourselves to the relatively
simple case of forms of weight two, and we study the differences
between the two $\ZZ$-structures, at primes whose square does not
divide the level. The main reason to study just this case is its
application in the article~\cite{AgSt1}, where our results are used to
obtain information about ``generalized Manin constants''. Another
reason is that more general cases are more difficult to treat, and
lead to results that are more complicated.

At the prime~$2$ the difference between the two $\ZZ$-structures that
we study can be bigger than at other primes. This special behaviour,
caused by the Hasse invariant being of weight one, is related to
spaces of modular forms modulo~$2$ of weight one (defined as by Katz,
i.e., mod~$2$ forms that cannot necessarily be lifted to
characteristic zero). As a byproduct of our work we find a description
of these spaces of weight one forms purely in terms of the integral
Hecke algebra associated to weight two cusp forms, together with the
Atkin-Lehner involution~$W_2$. These data can be computed from the
singular homology group $\rH_1(X_0(N)(\CC),\ZZ)$, for example using
modular symbols algorithms. Section~\ref{sec3} gives this description,
as an $\FF_2$-vector space equipped with Hecke operators $T_n$ for
odd~$n$.

In Section~\ref{sec4}, we give another (more general) description of
spaces of mod~$p$ modular forms of weight one in terms of a Hecke
algebra associated to forms of weight~$p$, based on properties of the
derivation~$qd/dq$. Theorem~\ref{thm4.7} gives a description of the
$\FF_p$-vector spaces of cusp forms of weight one on~$\Gamma_1(N)$
with a fixed character, equipped with all Hecke operators, in terms of
data that can be computed using modular symbols. Section~\ref{sec4a}
contains a result that says that some mod~$p$ Hecke algebras can be
computed using parabolic group cohomology with coefficients
mod~$p$. Finally, in Section~\ref{sec5}, we give a result about
eigenspaces.

It should be said that there are published methods about how one can
compute spaces of weight one modular forms, and even some tables. For
example there are the methods described in~\cite{Kiming1} (see also
\cite{BuzzardStein1}). But those methods are aimed at computing spaces
$S_1(\Gamma_1(N),\eps,\CC)$ of characteristic zero forms of weight
one. It seems that the methods as in~\cite{Kiming1} can be adapted to
compute spaces of Katz's modular forms mod~$p$, but that, at least for
small~$p$, they are less directly related to Hecke algebras that one
can easily compute using group cohomology or modular symbols
algorithms. 

The results in this article (especially the part about computing
spaces of mod~$p$ weight one forms) are not so original. The methods
employed are due for a large part to Katz, and Deligne and Rapoport,
and some ideas are already implicitly contained, in special cases,
in~\cite{BuzzardStein1}. But on the other hand there seem to be no
tables of mod~$p$ modular forms of weight one, and worse, no published
algorithm, to compute such tables. We hope that the results of the
last three sections will be used for computations. It would be
interesting to know if the examples of non-Gorenstein Hecke algebras
found by Kilford (see \cite{Kilford1} and~\cite[\S3.7.1]{RibetStein1})
can be ``explained'' by weight one phenomena. Another example is the
question whether the mod~$2$ Galois representation associated to a
mod~$2$ eigenform of weight one is unramified at~$2$, even if one is
in the exceptional case (see~\cite[Prop.~2.7]{Edix4}, the end of the
introduction of~\cite{Edix4}
and~\cite[Cor.~0.2]{ColemanVoloch}). Finally, one could test the
question of existence of companion forms in the exceptional case when
$p=2$ (see \cite{Gross1} and~\cite{ColemanVoloch}).

Already in 1987, Mestre wrote a letter to Serre in which he described
the results of computations that he had done on weight one forms
mod~$2$ on $\Gamma_0(N)$ with $N$ prime. His method (suggested to him
by Serre) is the one described here in detail (and with justification)
in Section~\ref{sec4} (although he probably computed the weight two
Hecke algebra using the so-called ``graph method''). The aim of
Mestre's computations was to test the idea that weight one forms
mod~$p$ as defined by Katz should correspond to odd irreducible Galois
representations mod~$p$ that are unramified at~$p$. The most
interesting results are some examples of mod~$2$ eigenforms of weight
one whose associated Galois representations have
image~$\SL_2(\FF_8)$. After this article (without the appendices) was
completed in March~2002, I asked Mestre about these examples. He then
sent me a copy of his letter, and it seemed appropriate to include it
as an appendix to this article. But then the results should be
verified in some independent way. That task has now been completed by
Gabor Wiese, using MAGMA and Stein's package HECKE. The second
appendix gives an account of Wiese's computations.

I thank Amod Agashe and William Stein for asking me about differences
between $\ZZ$-structures, and for motivating me enough to write up the
results below (and I apologize to them for taking so much time). I
thank Kevin Buzzard and Christophe Breuil for their corrections and
comments on a first version of this text.  I thank Gabor Wiese for
pointing out some mistakes, for discussions on this subject, and for
his work on the appendix that he has written. Finally, I thank
Jean-Fran{\c c}ois Mestre for letting me include his letter to Serre
as an appendix.

\section{Comparison of two integral structures}\label{sec2}
\subsectionrunin{}
Let $N\geq 1$ be an integer. For $A$ a subring of~$\CC$, we let
$S_2(A)=S_2(\Gamma_0(N),A)$ be the $A$-module of $f$ in
$S_2(\Gamma_0(N),\CC)$ whose $q$-expansion is in~$A[[q]]$. We want to
compare the $\ZZ$-module $S_2(\ZZ)$ with another one that comes from
algebraic geometry over~$\ZZ$. Let $X=X_0(N)$ be the modular curve
over~$\ZZ$ that is the compactified coarse moduli scheme for elliptic
curves with a given cyclic subgroup scheme of rank~$N$, as constructed
and described in the book by Katz and Mazur~\cite{KaMa}. Let $J$ be
the N\'eron model over~$\ZZ$ of the jacobian variety of~$X_\QQ$ (see
\cite{BLR1} for generalities about N\'eron models, and
\cite[Thm.~9.7]{BLR1} for some special results in the case of modular
curves).  Then $J$ is a smooth group scheme over~$\ZZ$, of relative
dimension the genus of~$X_\QQ$, and its cotangent space at the origin
$\Cot_0(J)$ is a $\ZZ$-structure on $\Cot_0(J_\QQ) =
\rH^0(X_\QQ,\Omega^1)$, with $\Omega^1$ the $\calO_{X_\QQ}$-module of
K\"ahler differentials. The Kodaira-Spencer isomorphism identifies
$S_2(\CC)$ with
$\rH^0(X_\CC,\Omega^1)=\CC\otimes\rH^0(X_\QQ,\Omega^1)$. We want to
compare the sub $\ZZ$-modules $S_2(\ZZ)$ and $\Cot_0(J)$
of~$S_2(\CC)$. In order to do this we use the following algebraic
interpretation of the $q$-expansion map.

The standard cusp~$\infty$, that corresponds to the generalized
elliptic curve $\PP^1_\ZZ$ with its points $0$~and~$\infty$ identified
and equipped with its subgroup scheme~$\mu_N$ in the terminology of
Deligne and Rapoport~\cite{DeRa}, is a $\ZZ$-valued point of~$X$,
i.e., an element of~$X(\ZZ)$. The Tate curve $\Gm{}{}/q^\ZZ$ over
$\ZZ((q))$ gives an isomorphism from the formal spectrum of $\ZZ[[q]]$
to the completion of~$X$ along~$\infty$ (\cite[Ch.~8]{KaMa} and
\cite[\S1.2]{Edix1}). In particular, the image of $\infty$ lies in the
open subscheme $X^\sm$ of $X$ on which the morphism to $\Spec(\ZZ)$ is
smooth. A differential form~$\omega$ on some open neighborhood
of~$\infty$ in~$X$ can be expanded (uniquely) as
$\sum_{n\geq1}a_n(\omega)q^n (dq/q)$, with the $a_n$ in~$\ZZ$. In the
same way, we have a $q$-expansion map from $\rH^0(X_\QQ,\Omega)$
to~$\QQ[[q]]$. 

The following result is well-known, but we give a proof anyway. It
will allow us to view $S_2(\ZZ)$ as the subset of elements of 
$\rH^0(X_\QQ,\Omega^1)$ that satisfy some integrality condition at all
prime numbers~$p$. 
\begin{subprop}\label{prop2.1}
The sub $\QQ$-vector spaces $S_2(\QQ)$ and $\rH^0(X_\QQ,\Omega^1)$ of
$S_2(\CC)$ are equal. 
\end{subprop}
\begin{proof}
We follow the arguments used by Katz in \cite[1.6]{Katz1}, and by
Deligne and Rapoport in \cite[VII, 3.9]{DeRa}. 

The $q$-expansion map from $\rH^0(X_\QQ,\Omega^1)$ to $\CC[[q]]$ has
image in~$S_2(\QQ)$. So it remains to prove the other inclusion. Let
$\sum a_nq^n$ be in $S_2(\QQ)$, and let $\omega=\sum a_nq^n\,dq/q$ in
$\Omega^1(X_\CC)$ be its associated one-form. The fact that
$\rH^0(X_\CC,\Omega^1)=\CC\otimes\rH^0(X_\QQ,\Omega^1)$ implies that
$\sum a_nq^n$ is actually in the subring $\CC\otimes\QQ[[q]]$
of~$\CC[[q]]$. We consider the following commutative diagram:
$$
\begin{CD}
0 @>>> \rH^0(X_\QQ,\Omega^1) @>>> \CC\otimes\rH^0(X_\QQ,\Omega^1) @>>>
(\CC/\QQ)\otimes\rH^0(X_\QQ,\Omega^1) @>>> 0 \\
@. @VVV @VVV @VVV @. \\
0 @>>> \QQ[[q]] @>>> \CC\otimes\QQ[[q]] @>>> (\CC/\QQ)\otimes\QQ[[q]]
@>>> 0 
\end{CD}
$$
The two rows are exact, and the three vertical arrows are injective
(use that $X_\QQ$ is integral, that $\Omega^1_{X_\QQ/\QQ}$ is a line
sheaf, and that the functors $\CC\otimes-$ and $(\CC/\QQ)\otimes-$ are
exact). The statement is now obvious.
\end{proof}

\begin{subprop}\label{prop2.2}
The sub $\ZZ$-module $\Cot_0(J)$ of $S_2(\QQ)$ is contained
in~$S_2(\ZZ)$. 
\end{subprop}
\begin{proof}
We have to show that the $q$-expansion of an element of~$\Cot_0(J)$,
viewed as an element of $\rH^0(X_\QQ,\Omega^1)$, lies
in~$\ZZ[[q]]$. Evaluation at $0$ identifies $\rH^0(J_\QQ,\Omega^1)$
with $\Cot_0(J_\QQ)$, and $\rH^0(J,\Omega^1)$ with $\Cot_0(J)$ (use
that the elements of $\rH^0(J_\QQ,\Omega^1)$ are translation
invariant). The identification of $\rH^0(J_\QQ,\Omega^1)$ with
$\rH^0(X_\QQ,\Omega^1)$ is given by pullback via the morphism
$X_\QQ\to J_\QQ$ that sends an $S$-valued point $P$ of $X_\QQ$ to the
class of the invertible $\calO_{X_S}$-module $\calO_{X_S}(P-\infty)$,
where $S$ is any $\QQ$-scheme. By the N\'eron property of~$J$, this
morphism extends (uniquely) to a morphism $X^\sm\to J$. It follows
that $\rH^0(J,\Omega^1)$ (and hence $\Cot_0(J)$) is mapped into
$\rH^0(X^\sm,\Omega^1)$, and hence, via the $q$-expansion map,
into~$\ZZ[[q]]$.
\end{proof}
\subsectionrunin{}
Our aim is to compare $S_2(\ZZ)$ and~$\Cot_0(J)$, so we have to
describe the quotient $S_2(\ZZ)/\Cot_0(J)$. As this quotient is
torsion, it is the direct sum over all prime numbers $p$ of its
$p$-primary part. We will give a description of the $p$-primary part
for $p$ such that $p^2$ does not divide~$N$ (recall that
$X=X_0(N)$). 

By definition, $S_2(\ZZ)$ is the subset of elements of
$\rH^0(X_\QQ,\Omega^1)$ that satisfy certain integrality conditions at
all prime numbers~$p$. More precisely, for $\omega$ in $S_2(\ZZ)$ and
$p$ prime, the condition is that the $q$-expansion of $\omega$ at the
cusp~$\infty$, which is an element $(\sum_{n\geq1}a_nq^n)dq/q$ of
$\QQ[[q]]{\cdot}dq$, satisfies $v_p(a_n)\geq0$ for all~$n$. Our first
step is to reinterpret these integrality conditions geometrically, as
in~\cite[VII, Thm.~3.10]{DeRa}.

Let $\eta$ denote the generic point of~$X$. We note that
$\Omega^1_{X^\sm/\ZZ}$ is an invertible $\calO_{X^\sm}$-module. Hence,
to each element $\omega$ of $\Omega^1_{X/\ZZ,\eta}$ we can associate
its multiplicity $v_C(\omega)$ along each prime divisor $C$
of~$X^\sm$. This multiplicity is defined as~$v_{\eta_C}(f)$, where
$\eta_C$ is the generic point of~$C$, $v_{\eta_C}$ the discrete
valutation on the discrete valuation ring~$\calO_{X,\eta_C}$, and
$\omega=f{\cdot}\omega_C$ with
$\Omega^1_{X^\sm/\ZZ,\eta_C}=\calO_{X,\eta_C}\omega_C$. In these
terms, we have:
$$
\rH^0(X^\sm,\Omega^1_{X^\sm/\ZZ}) = \{\omega\in
\rH^0(X_\QQ,\Omega^1_{X_\QQ/\QQ})\;|\;
\text{$v_C(\omega)\geq0$ for all $C$.}\}
$$
In order to relate the integrality conditions on $q$-expansions to the
$v_C(\omega)$, we look at what happens under completion at some closed
point.

Let $x$ be a closed point of~$X^\sm$. Then $\calO_{X,x}$ is a
two-dimensional regular noetherian local ring, and since
$X\to\Spec(\ZZ)$ is smooth at~$x$, it has a system of parameters of
the form $(p,t)$, with $p$ a prime number. The completion
$\calO_{X,x}^\wedge$ is then isomorphic to~$W(k)[[t]]$, with $W(k)$
the ring of Witt vectors of the residue field $k$ at~$x$. Let $C$ be
the irreducible component of $X_{\FF_p}$ that $x$ lies on. Then
$\calO_{C,x}=\calO_{X,x}/p\calO_{X,x}$ and the morphism from
$\calO_{C,x}$ to its completion $\calO_{C,x}^\wedge=k[[t]]$ is
injective. Let now $f$ be in~$\calO_{X,x}$, and put $m:=v_C(f)$. Then,
as $p$ is a prime element of $\calO_{X,x}$ and a uniformizer
of~$\calO_{X,\eta_C}$, we have $f=p^mf'$, with $f'$ in $\calO_{X,x}$
such that the image $\ol{f'}$ of $f'$ in $\calO_{C,x}$ is not
zero. Let us write the image of $f$ in $W(k)[[t]]$ as $\sum a_n t^n$
(with the $a_n$ in~$W(k)$). Then we see that $m=\min_n v_p(a_n)$ (use
that the image of $\ol{f'}$ in $k[[t]]$ is non-zero).

If we take for $x$ the cusp $\infty$ in $X(\FF_p)$ for some prime
number~$p$, the discussion above gives that the integrality condition
at $p$ for $\omega$ in $\rH^0(X_\QQ,\Omega^1_{X_\QQ/\QQ})$ to be in
$S_2(\ZZ)$ is just $v_C(\omega)\geq0$, with $C$ the irreducible
component of $X_{\FF_p}$ that contains~$\infty$. This proves the
following proposition. 

\begin{subprop}\label{prop2.3}
Let $\omega$ be in $\rH^0(X_\QQ,\Omega^1_{X_\QQ/\QQ})$. Then $\omega$
is in $S_2(\ZZ)$ if and only if for every prime number~$p$ it has
multiplicity $\geq 0$ along the irreducible component of $X_{\FF_p}$
that contains~$\infty$. In other words:
$$
S_2(\ZZ)=\rH^0(X_\infty,\Omega^1_{X/\ZZ}), 
$$
with $X_\infty$ the complement in $X$ of the union of all the
irreducible components of all $X_{\FF_p}$ that do \emph{not}
contain~$\infty$. 
\end{subprop}
\subsectionrunin{}\label{sec2.6}
It may come as a surprise that, with this characterisation, $S_2(\ZZ)$
is a finitely generated $\ZZ$-module, since $X_\infty$ is not proper
if~$N\neq1$. But it is not hard to show that for $C$ and $C'$
irreducible components of the same~$X_{\FF_p}$, the differences
$|v_C(\omega)-v_{C'}(\omega)|$, for $\omega\neq0$ in
$\rH^0(X_\QQ,\Omega^1_{X_\QQ/\QQ})$, are bounded uniformly
in~$\omega$. In fact, after choosing an extension of
$\Omega^1_{X_\infty/\ZZ}$ to an invertible
$\calO_X$-module~$\mathcal{L}$ (after a resolution of singularities,
if necessary), one obtains a bound that depends only on the
combinatorial data: the dual graph plus intersection numbers,
multiplicities, genera and degrees of restrictions of $\mathcal{L}$
associated to $X_{\FF_p}$ and~$\mathcal{L}$. Such bounds imply that
$S_2(\ZZ)$ is a lattice in $\rH^0(X_\QQ,\Omega^1_{X_\QQ/\QQ})$. The
computations that we will do below are an explicit and exact version
of the proof just alluded to, but in a simple case. We note that in
\cite[VII, \S3]{DeRa} similar computations have been done, but for
modular forms of arbitrary weight~$k$, seen as sections
of~$\underline{\omega}^{\otimes k}$.

If $p$ is a prime that does not divide~$N$, then $X_{\FF_p}$ is
irreducible, and hence the quotient $S_2(\ZZ)/\Cot_0(J)$ is trivial
at~$p$. If $p$ is a prime number such that $p^2$ divides $N$ we do not
want to say anything about the $p$-part of the quotient in this
article. A good reason for that is that the results become much more
difficult to describe. One finds some computations in
\cite[\S4.6]{Edix2} and \cite[\S4]{Edix3} for the case where $p^2$
exactly divides~$N$.

So suppose from now on that $p$ is a prime dividing $N$ exactly, i.e.,
such that $N/p$ is not divisible by~$p$. Then $X$ is semistable
at~$p$; in particular, $X_{\ZZ_p}$ is normal. Its fibre $X_{\FF_p}$
has two irreducible components, $C_\infty$~and~$C_0$, both isomorphic
to~$X_0(N/p)_{\FF_p}$ (and hence smooth), with $C_\infty$ containing
the cusp~$\infty$, and hence $C_0$ the cusp~$0$. The intersection
$C_\infty\cap C_0$ consists of the supersingular points on $C_\infty$
and~$C_0$, and the intersections are transversal. We identify
$C_\infty$ with $X_0(N/p)_{\FF_p}$ via the following
construction: 
$$
X_0(N/p)_{\FF_p}\lto C_\infty,\quad (E/S/\FF_p,G)\mapsto
(E/S/\FF_p,G,\ker{F}), 
$$
where $E/S/\FF_p$ is an elliptic curve over an $\FF_p$-scheme, where
$G$ is a cyclic subgroup scheme of rank $N/p$ of~$E/S$, and where
$F\colon E\to E^{(p)}$ is the Frobenius morphism of $E$
over~$S$. Likewise, we identify $C_0$ with $X_0(N/p)_{\FF_p}$ as
follows: 
$$
X_0(N/p)_{\FF_p}\lto C_\infty,\quad (E/S/\FF_p,G)\mapsto
(E^{(p)}/S/\FF_p,G^{(p)},\ker{V}), 
$$
where $V\colon E^{(p)}\to E$ is the Verschiebung (the dual of~$F$). 

Let $\Omega$ be the dualising sheaf for $X_{\ZZ_p}$ over~$\ZZ_p$ (see
\cite{MazurRibet1} for a discussion of this sheaf in the context of
semi-stable modular curves). Then $\Omega$ is the direct image of
$\Omega^1_{X^\sm_{\ZZ_p}/\ZZ_p}$ via the inclusion of $X^\sm_{\ZZ_p}$
into $X_{\ZZ_p}$, and it is an invertible
$\calO_{X_{\ZZ_p}}$-module. We have:
$$
\Cot_0(J_{\ZZ_p}) = \rH^1(X_{\ZZ_p},\calO)^\vee = \rH^0(X_{\ZZ_p},\Omega). 
$$ 
Let $\omega$ be a non-zero element of $\rH^0(X_{\QQ_p},\Omega)$. We
want to bound $|v_{C_\infty}(\omega)-v_{C_0}(\omega)|$, since the
exponent of $p$ in the exponent of the quotient $S_2(\ZZ)/\Cot_0(J)$
is the maximal value of this expression over all~$\omega$. The
automorphism $W_N$ of $X$ interchanges $\infty$ and~$0$, so we may as
well suppose that $v_{C_0}(\omega)=0$ and that
$m:=v_{C_\infty}(\omega)\geq0$. Then the restriction $\omega|_{C_0}$
is a non-zero section of~$\Omega|_{C_0}$, and the fact that $\omega$
has a zero of order~$m$ along~$C_\infty$ implies that $\Omega|_{C_0}$
has zeros of at least some order at the supersingular points. As the
number of zeros cannot exceed the degree of~$\Omega|_{C_0}$, we will
get an upper bound for~$m$.

There is now a minor problem: $X_{\ZZ_p}$ is not necessarily regular,
and also, the degree of $\Omega|_{C_0}$ is easier to relate to the
number of supersingular points if we are working on a fine moduli
space. So we choose to work on some cover $X'$ of~$X_{\ZZ_p}$, say
obtained by adding, to the $\Gamma_0(N)$-structure, a
$\Gamma(3)$-structure when $p\neq3$ or a $\Gamma(4)$-structure if
$p=3$.  Note that $X'\to X_{\ZZ_p}$ is finite, and etale outside the
locus of $j$-invariants~$0$, $1728$~and~$\infty$. Let $C_\infty'$ and
$C_0'$ be the inverse images of $C_\infty$ and $C_0$, respectively,
with their reduced scheme structure. Then the pullback of $\omega$ to
$X'$ has valuation zero along every irreducible component of~$C_0'$,
and a zero of order~$m$ along every irreducible component
of~$C_\infty'$. As $X'$ is regular, this implies that $\omega|_{C_0'}$
is a non-zero global section of~$\Omega|_{C_0'}$, with zeros of order
at least $m$ at all supsersingular points. In other words, $\omega$ is
a non-zero global section of $\Omega|_{C_0'}(-mS)$, where $S$ is the
divisor given by the sum of all supersingular points. Hence
$\deg(\Omega|_{C_0'}(-mS))\geq0$. Let ``$\cusps$'' denote the sum of
the cups on~$C_0'$. Then we have:
\begin{align*}
\deg(\Omega|_{C_0'}(-mS)) & = \deg(\Omega^1_{C_0'/\FF_p}) + \# S -
m{\cdot}\# S \\
& = 2{\cdot}\deg(\underline{\omega}) - \#\cusps + (1-m)\# S \\
& = \frac{2}{p-1}\# S - \#\cusps + (1-m)\# S, 
\end{align*}
where we have used the equality
$\Omega|_{C_0'}=\Omega^1_{C_0'/\FF_p}(S)$, the Kodaira-Spencer
isomorphism
$\Omega^1_{C_0'/\FF_p}(\cusps)\to\underline{\omega}^{\otimes 2}$ (for
the family of elliptic curves obtained via the isomorphism
with~$C_\infty'$; see below for an explanation), and the Hasse
invariant in $\rH^0(C_0',\underline{\omega}^{\otimes p-1})$, whose
divisor is~$S$. It follows that:
$$
m < 1 + \frac{2}{p-1}, 
$$
and we have proved the following result. 
\begin{subprop}\label{prop2.4}
Let $N\geq 1$ be an integer and let $p$ be a prime that exactly
divides~$N$. Then the $p$-part of the quotient $S_2(\ZZ)/\Cot_0(J)$ is
annihilated by~$p$, if $p>2$, and by $4$ if $p=2$.
\end{subprop}
\subsectionrunin{}
Our next aim is to describe the cokernel more precisely, at least if
$p\neq 2$. But let us first compare Proposition~\ref{prop2.4} with
\cite[VII, Prop.~3.20]{DeRa}. It is clear that that result, stated for
modular forms on~$\Gamma_0(p)$, remains true if one adds prime to $p$
level structure (the same proof works). It says that for a non-zero
weight $k$ form on~$\Gamma_0(N)$ with coefficients in~$\QQ$,
corresponding to an element $f$ of
$\rH^0(X',\underline{\omega}^{\otimes k})$, one has:
$$
\left|v_{C_\infty'}(f) - v_{C_0'}(f) + \frac{k}{2}\right| \leq
 \frac{1}{2} {\cdot} k {\cdot} \frac{p+1}{p-1}. 
$$ 
This result is not symmetric in $C_\infty$ and $C_0$ because the sheaf
$\underline{\omega}$ is not: the $j$-invariant is separable
on~$C_\infty$, and inseparable on~$C_0$. Proposition~\ref{prop2.4} can
be deduced from this inequality, if one takes into account in its
proof that $f$ is a cuspform, and that the Kodaira-Spencer morphism
$\underline{\omega}^{\otimes 2}\to\Omega^1(\cusps)$ on $X'$ has a zero
of order one along~$C_0'$. This last fact can be seen by looking at
the cusps, via the Tate curve:
$$
(dt/t)^{\otimes 2}\mapsto dq/q = p{\cdot}d(q^{1/p})/q^{1/p}.
$$

Let us now try to determine the $p$-part of $S_2(\ZZ)/\Cot_0(J)$ more
precisely. Both $S_2(\ZZ)$ and $\Cot_0(J)$ are stable under all 
endomorphisms ~$T_n$, $n\geq1$, of $S_2(\Gamma_0(N),\CC)$. For
$S_2(\ZZ)$ this follows from the usual formulas in terms of
$q$-expansions (\cite[(12.4.1)]{DiIm}): 
$$
a_m(T_nf)=\sum_{\substack{0<d|(n,m)\\(d,N)=1}} d^{k-1}a_{nm/d^2}(\ld d\rd f), 
$$
for $f$ in $S_k(\Gamma_1(N),\CC)$, $n$ and $m$ positive integers. By
construction, $\Cot_0(J)$ is stable under all endormorphisms
of~$J_\QQ$, hence in particular under the~$T_n$. The Atkin-Lehner
involutions~$W_d$, for $d$ dividing $N$ such that $d$ and $N/d$ are
relatively prime, stabilize $\Cot_0(J)$, but do not necessarily
stabilize $S_2(\ZZ)$, since they may interchange $C_\infty$
and~$C_0$. In fact, $W_d$ interchanges $C_\infty$ and $C_0$ if and
only if $p$ divides~$d$. 

By Proposition~\ref{prop2.4}, the cokernel of:
$$
\rH^0(X_{\ZZ_p},\Omega) \longrightarrow 
\{\omega \in \rH^0(X_{\QQ_p},\Omega)\;|\;v_{C_0}(\omega)\geq -1\}
$$
is the $p$-part of $S_2(\ZZ)/\Cot_0(J)$ if $p\neq 2$, and the piece
killed by $p$ of it if $p=2$. If $\omega$ is in
$\rH^0(X_{\QQ_p},\Omega)$ such that $v_{C_\infty}(\omega)\geq0$ and
$v_{C_0}(\omega)\geq -1$, then $p\omega |_{C_0}$ is an element of
$\rH^0(C_0,\Omega^1_{C_0/\FF_p})$. We have an exact sequence:
$$
0\to \rH^0(X_{\ZZ_p},\Omega) \to
\{\omega \in \rH^0(X_{\QQ_p},\Omega)\;|\;v_{C_0}(\omega)\geq -1\} 
\to \rH^0(C_0,\Omega^1_{C_0/\FF_p}).
$$
On the other hand, suppose $\omega$ is in
$\rH^0(C_0,\Omega^1_{C_0/\FF_p})$. Then we have an element, call
it~$(0,\omega)$, of $\rH^0(X_{\FF_p},\Omega)$ whose restrictions
to~$C_\infty$ and~$C_0$ are $0$~and~$\omega$, respectively. The fact
that $\Omega$ is the dualising sheaf implies that the formation of its
cohomology commutes with the base change $\ZZ_p\to\FF_p$, and hence
that $(0,\omega)$ is the reduction mod~$p$ of an element,
$\widetilde{\omega}$~say, in~$\rH^0(X_{\ZZ_p},\Omega)$. Then $\omega$
is the image of $p^{-1}\widetilde{\omega}$, in the exact sequence
above. We have now proved the first claim of the following
proposition.

\begin{subprop}\label{prop2.5}
Let $N\geq1$ be an integer and let $p$ be prime that exactly
divides~$N$. Then the construction above gives an isomorphism between
$\rH^0(C_0,\Omega^1_{C_0/\FF_p})$ and the $p$-part of
$S_2(\ZZ)/\Cot_0(J)$ if $p>2$, and its part killed by $2$ if $p=2$.
We identify $\rH^0(C_0,\Omega^1_{C_0/\FF_p})$ and $\FF_p\otimes
S_2(\Gamma_0(N/p),\ZZ)$, using the chosen isomorphism between
$X_0(N/p)_{\FF_p}$ and~$C_0$. Then the endomorphism $T_n$ of
$S_2(\ZZ)/\Cot_0(J)$ induces $T_n$ on $\FF_p\otimes
S_2(\Gamma_0(N/p),\ZZ)$ if $n$ is prime to~$p$, and $0$ if $n$ is a
multiple of~$p$.
\end{subprop}
\begin{proof}
It remains to prove the second claim. Let $l$ be a prime number, and
let $\omega$ be in $\rH^0(C_0,\Omega^1_{C_0/\FF_p})$. Going through
the construction of the isomorphism above, one finds that the
endomorphism $T_l$ of $S_2(\ZZ)$ sends $\omega$ to
$(T_l(0,\omega))|_{C_0}$, where $(0,\omega)$ is the element of
$\rH^0(X_{\FF_p,\Omega})$ with restrictions $0$ and~$\omega$ to
$C_\infty$ and~$C_0$, respectively. If $l$ is not equal to~$p$, this
means that $T_l$ induces the usual $T_l$ on $\FF_p\otimes
S_2(\Gamma_0(N/p),\ZZ)$. For $l=p$, one notes (see
\cite[\S6.6]{Edix4}) that $T_p$ induces the Frobenius endomorphism of
the jacobian of~$C_0$. 
\end{proof}
Finally, we investigate what happens at~$2$. We will prove the
following result.

\begin{subprop}\label{prop2.6}
Let $N\geq 1$ be an integer. Suppose that $2$ exactly divides~$N$, and
that $N$ is divisible by a prime number $q\equiv-1$ modulo~$4$ (this
last condition means that no elliptic curve over~$\Fbar_2$ equipped
with a cyclic subgroup of order~$N/2$ has an automorphism of
order~$4$). Let $M$ denote the $2$-part of~$S_2(\ZZ)/\Cot_0(J)$. Then:
\begin{align*}
& M[4]  = M,\\
& M/M[2] = \rH^0(C_0,\Omega^1_{C_0/\FF_2}(-S)) = 
S_1(\Gamma_0(N/2),\FF_2)_\Katz,\\
& M[2] = \rH^0(C_0,\Omega^1_{C_0/\FF_2}),
\end{align*}
where $S_1(\Gamma_0(N/2),\FF_2)_\Katz$ denotes the $\FF_2$-vector
space of weight one cusp forms on $\Gamma_0(N/2)$ over~$\FF_2$,
defined as by Deligne and Katz (see \cite{Katz1}, \cite[\S2]{Edix4} or
\cite[\S1]{Edix5} and the references therein for a definition and
properties of Katz's modular forms). For $n\geq1$, the endomorphism
$T_n$ of $S_2(\ZZ)$ induces the endomorphism $T_n$ of
$S_1(\Gamma_0(N/2),\FF_2)_\Katz$ if $n$ is odd, and $0$ if $n$ is
even.
\end{subprop}
\begin{proof}
The statements concerning $M[4]$ and $M[2]$ have already been proved
in Propositions~\ref{prop2.4} and~\ref{prop2.5}, and hence are even
valid without the extra condition on~$N/2$. It remains to prove the
description of~$M/M[2]$, and the action of $T_n$ on it.

Suppose that $\omega$ is an element of $\rH^0(X_{\QQ_2},\Omega)$ with
$v_{C_\infty}(\omega)\geq0$. Then $v_{C_0}(\omega)\geq -2$, hence
$4\omega$ is an element of $\rH^0(X_{\ZZ_2},\Omega)$ with at least a
double zero along~$C_\infty$. A local computation shows that
$4\omega|_{C_0}$ is an element of
$\rH^0(C_0,\Omega^1_{C_0/\FF_2}(-S))$. We get a map:
\begin{equation}\label{eqn2.6.1}
\{\omega \in \rH^0(X_{\QQ_2},\Omega)\;|\;v_{C_0}(\omega)\geq -2\} 
\lto \rH^0(C_0,\Omega^1_{C_0/\FF_2}(-S)), 
\end{equation}
with kernel the submodule of $\omega$ with $v_{C_0}(\omega)\geq -1$.
We will show that this map is surjective. For this we will work on a
suitable cover of~$X_{\ZZ_2}$.  Let $\ZZ_4$ be the unramified
quadratic extension of~$\ZZ_2$. Then the map~\ref{eqn2.6.1} is
surjective if and only if its analog after base change to~$\ZZ_4$
is. Fix an element $\zeta$ of order three in~$\ZZ_4^*$ and let $X''\to
X_{\ZZ_4}$ be the cover obtained by adding a full level three
structure with Weil pairing~$\zeta$. Let $X'\to X_{\ZZ_4}$ be the
cover obtained by dividing out the action of the Sylow $2$-group
of~$\SL_2(\FF_3)$. Then $X_{\ZZ_4}$ is the quotient of $X'$ by the
action of the quotient $G$ of order three of~$\SL_2(\FF_3)$. Our
hypothesis that $N$ is divisible by a prime number that is $-1$
modulo~$4$ implies that $X''\to X'$ is etale, hence that $X'$ is
regular. We denote by $C_0'$ and $C_\infty'$ the inverse images of
$C_0$ and~$C_\infty$, with their reduced scheme structure. Then $C_0'$
and $C_\infty'$ are smooth geometrically irreducible curves.

Let $\omega$ be in $\rH^0(C_0,\Omega^1_{C_0/\FF_2}(-S))$. The inverse
image of $\omega$ in $\rH^0(C_0',\Omega^1_{C_0'/\FF_4}(-S))$ will
still be denoted by~$\omega$. Suppose that we have an element
$\widetilde{\omega}$ of $\rH^0(X',\Omega)^G$. Then there is a unique
$\eta$ in $\rH^0(U,\Omega)$, with $U$ the open subset of $X_{\ZZ_4}$
over which $\pi\colon X'\to X_{\ZZ_4}$ is etale, such that
$\widetilde{\omega}=\pi^*\eta$. Normality of $X_{\ZZ_4}$ and the fact
that $\eta$ extends in codimension one imply that $\eta$ is in
$\rH^1(X_{\ZZ_4},\Omega)$ (we note that we did not use that the
$G$-action is tame). Hence it suffices to show that there is a
$G$-invariant element $\widetilde{\omega}$ of $\rH^0(X',\Omega)$ that
has at least a double zero along~$C_\infty'$, and whose restriction
to~$C_0'$ is~$\omega$.  The last two of these properties of such an
$\widetilde{\omega}$ can already be seen from its restriction
to~$X'_{\ZZ_4/4\ZZ_4}$, and if $\widetilde{\omega}$ does satisfy them,
then so does its projection to the $G$-invariants (here we use that
$G$ has order prime to~$2$).  As $\rH^1(X',\Omega)$ is torsion free,
the reduction map from $\rH^0(X',\Omega)$ to
$\rH^0(X'_{\ZZ_4/4\ZZ_4},\Omega)$ is surjective. Let $(0,\omega)$ be
the section of $\Omega$ on the closed subscheme of $X'$ defined
by~$I_{C_\infty'}^2I_{C_0'}$ that coincides with $\omega$
modulo~$I_{C_0'}$ and is zero modulo~$I_{C_\infty'}^2$. It is now
enough to show that $(0,\omega)$ can be lifted to a section
of~$\Omega$ over~$X'_{\ZZ_4/4\ZZ_4}$. The short exact sequence of
sheaves corresponding to our lifting problem is:
$$
0\to 2I_{C_0'}/4\calO_{X'} \to \calO_{X'}/4\calO_{X'} \to 
\calO_{X'}/2I_{C_0'}\to 0, 
$$
tensored over $\calO_{X'}$ with~$\Omega$. So all we have to show is
that the map from $\rH^1(X',\Omega\otimes 2I_{C_0'}/4\calO_{X'})$ to
$\rH^1(X'_{\ZZ_4/4\ZZ_4},\Omega)$ is injective. One verifies that
$\Omega\otimes 2I_{C_0'}/4\calO_{X'} = \Omega^1_{C_0'/\FF_4}$, so that
its $\rH^1$ is~$\FF_4$, and that
$\rH^1(X',\Omega\otimes\calO_{X'}/2I_{C_0'})$ also has dimension one
over~$\FF_4$ (use a suitable short exact sequence with terms
$\calO_{X'}/I_{C_0'}$, $\calO_{X'}/2I_{C_0'}$
and~$\calO_{X'}/2\calO_{X'}$).  We have now proved that $M/M[2]$ is
the same as~$\rH^0(C_0,\Omega^1_{C_0/\FF_2}(-S))$.

Let us now prove the equality between
$\rH^0(C_0,\Omega^1_{C_0/\FF_2}(-S))$
and~$S_1(\Gamma_0(N/2),\FF_2)_\Katz$. Let $C_0''\to C_0$ be the cover
obtained by adding a full level three structure
to~$\Gamma_0(N/2)$. Then, by our hypotheses on~$N$, the morphism
$C_0''\to C_0$ is not wildly ramified. Let now $G$ denote the
group~$\GL_2(\FF_3)$, which acts on~$C_0''$, and on the universal
generalized elliptic curve with level structure over~it. It follows
that:
$$
\rH^0(C_0,\Omega^1_{C_0/\FF_2}(-S)) = 
\rH^0(C_0'',\Omega^1_{C_0''/\FF_2}(-S''))^G.
$$
By definition, we have: 
$$
S_1(\Gamma_0(N/2),\FF_2)_\Katz =
\rH^0(C_0'',\underline{\omega}(-\cusps))^G.
$$
The Kodaira--Spencer isomorphism and the Hasse invariant, which both
exist on~$C_0''$ and are $G$-equivariant, finish the proof of the
equality.

The determination of the endomorphism of
$\rH^0(C_0,\Omega^1_{C_0/\FF_2}(-S))$ induced by $T_l$ (with $l$ a
prime number) is done in the same way as in the proof
of~\ref{prop2.5}, and we don't repeat it. To describe it on
$S_1(\Gamma_0(N/2),\FF_2)_\Katz$ one uses the formulas relating Hecke
operators and $q$-expansions, and one uses that $T_2$ acts as zero. 
\end{proof}

\section{Computation of spaces of forms of weight one: first method.}
\label{sec3}
\subsectionrunin{} 
Proposition~\ref{prop2.6} leads to an algorithm to compute, for each
$N$ satisfying the hypotheses of that proposition, the $\FF_2$-vector
space $S_1(\Gamma_0(N/2),\FF_2)_\Katz$ with the action of $T_n$ for
all $n$ prime to~$2$. Let us explain this.

Let $N$ satisfy the hypotheses of Proposition~\ref{prop2.6}, and let
$\TT$ be the Hecke algebra associated to $S_2(\Gamma_0(N),\ZZ)$, i.e.,
the sub $\ZZ$-algebra of $\End(S_2(\Gamma_0(N),\ZZ))$ (or of
$\End_\CC(S_2(\Gamma_0(N),\CC))$, it gives the same algebra) generated
by all~$T_n$, $n\geq1$. Then we have a $\ZZ$-valued pairing between
$S_2(\Gamma_0(N),\ZZ)$ and~$\TT$, given by $(f,t)\mapsto a_1(tf)$. As
$a_1(T_nf)=a_n(f)$ for all $n$ and~$f$, this pairing is
perfect. Hence, as a $\TT$-module, $S_2(\Gamma_0(N),\ZZ)$ is
isomorphic to the $\ZZ$-dual~$\TT^\vee$ of~$\TT$. Let $W_2$ be the
Atkin-Lehner involution that acts on $S_2(\Gamma_0(N),\CC)$; as it is
induced by an automorphism of $X_0(N)_{\ZZ_2}$ that interchanges the
components $C_\infty$ and~$C_0$, we have:
$$
W_2 (\ZZ_2\otimes S_2(\Gamma_0(N),\ZZ)) =
\{\omega\in\rH^0(X_0(N)_{\QQ_2},\Omega^1)\;|\;v_{C_0}(\omega)\geq0\}.
$$
Hence: 
$$
\ZZ_2\otimes(S_2(\Gamma_0(N),\ZZ)\cap W_2 S_2(\Gamma_0(N),\ZZ)) = 
\rH^0(X_0(N)_{\ZZ_2},\Omega) = \ZZ_2\otimes\Cot_0(J).
$$
We define a $\TT$-module $M$ by:
\begin{equation}
M:= S_2(\Gamma_0(N),\ZZ)/(S_2(\Gamma_0(N),\ZZ)\cap W_2
S_2(\Gamma_0(N),\ZZ)) = \TT^\vee/(\TT^\vee\cap W_2\TT^\vee),  
\end{equation}
where the last intersection can be taken in~$\QQ\otimes\TT^\vee$.
Then, by Proposition~\ref{prop2.6} we have an isomorphism of
$\FF_2$-vector spaces:
\begin{equation}
S_1(\Gamma_0(N/2),\FF_2)_\Katz = M[4]/M[2], 
\end{equation}
such that $T_n$ on $M[4]/M[2]$ induces $T_n$ on
$S_1(\Gamma_0(N/2),\FF_2)_\Katz$ if $n$ is odd, and $0$ if $n$ is
even. 

\subsectionrunin{}
The theory of modular symbols (see \cite{Cremona1}, \cite{Merel1} or
William Stein's modular forms database web pages) allows one to
compute with $\rH_1(X_0(N)(\CC),\ZZ)$. In particular, one can compute
the Hecke algebra~$\TT$ as well as the involution~$W_2$. Here it is
useful to note that the $T_n$ with $n\leq
6^{-1}N\prod_{p|N}(1+p^{-1})$ generate $\TT$ as a $\ZZ$-module:
see~\cite{AgSt2} (which is a simple application
of~\cite{Sturm}). Hence one can compute $M$ and so one gets
$S_1(\Gamma_0(N/2),\FF_2)_\Katz$ with the Hecke operators $T_n$ with
$n$ odd. 

The methods of this section can be generalized to mod $2$ modular
forms on $\Gamma_1(N)$, with $N$ odd, by considering the space
$S_2(\Gamma_1(2N),\CC)$. They can also be generalized to an arbitrary
prime number~$p$, by considering the spaces
$S_p(\Gamma_0(p)\cap\Gamma_1(N),\CC)$ with $N$ prime to~$p$. But
because of the next section, we think that this is not worth the
effort. The results of this section can be used to check some special
cases of the results of the next section.

In the next section we will present yet another way to compute more
generally the spaces $S_1(\Gamma_1(N),\eps,\FF_p)_\Katz$, this time
with all Hecke operators, and also in a way that should allow for a
simple implementation starting from William Stein's Magma package
Hecke.

In a sense, our method described in this section uses the Hasse
invariant, which is a very natural modular form modulo~$p$. As it is
of weight $p-1$, multiplication by it gives a way for passing from
weight one forms to weight $p$ forms. In the next section, we also use
that if $f$ is a weight one form, then $f^p$ is a weight $p$ form.

\section{Spaces of forms of weight one: second method.}
\label{sec4}
\subsectionrunin{} Let $p$ be a prime number, and $N\geq1$ an integer
prime to~$p$. For any integer $k$ and any finite extension $\FF$ of
$\FF_p$ we have the $\FF$-vector space $S_k(\Gamma_1(N),\FF)_\Katz$ of
cuspidal modular forms on $\Gamma_1(N)$ of weight~$k$ and with
coefficients in~$\FF$; we refer (again) to \cite[\S2]{Edix4} or
\cite[\S1]{Edix5} for a definition and a description. These spaces are
associated to the moduli problem~$[\Gamma_1(N)]_\FF$ of elliptic
curves with a given point of order~$N$. Because of rationality
properties of the unramified cusps, we will also consider the spaces
$S_k(\Gamma_1(N),\FF)'_\Katz$ of Katz modular forms associated to the
moduli problem $[\Gamma_1(N)]'_\FF$ of elliptic curves with an
embedding of the group scheme~$\mu_N$ as for example
in~\cite{Gross1}. It is an advantage to have the unramified cusps
rational, because the Hecke correspondences $T_n$ decrease the
ramification of the cusps, hence map unramified cusps to unramified
cusps. We let $X_1(N)'_{\ZZ[1/N]}$ be the compactified modular curve
associated to $[\Gamma_1(N)]'_{\ZZ[1/N]}$. As the group schemes
$(\ZZ/N\ZZ)_\ZZ$ and $\mu_N$ become canonically isomorphic over
$\ZZ[1/N,\zeta_N]$, the same holds for the two stacks
$[\Gamma_1(N)]_{\ZZ[1/N]}$ and $[\Gamma_1(N)]'_{\ZZ[1/N]}$. It follows
that $S_k(\Gamma_1(N),\Fbar_p)_\Katz$ and
$S_k(\Gamma_1(N),\Fbar_p)'_\Katz$ are isomorphic (via the choice of an
element of order $N$ in~$\Fbar_p^*$), compatibly with the action of
Hecke, diamond and Atkin-Lehner operators. In particular, the Hecke
operators~$T_n$ ($n\geq1$) and diamond operators $\ld a\rd$ ($a$ in
$(\ZZ/N\ZZ)^*$) generate the same algebra
over~$\FF_p$. Proposition~\ref{prop4.4} shows that, for any
$\ZZ[1/N]$-algebra $R$, the Hecke modules $S_k(\Gamma_1(N),R)_\Katz$
and $S_k(\Gamma_1(N),R)'_\Katz$ are isomorphic, but we will not use
this fact. We have:
$$
S_k(\Gamma_1(N),\FF)'_\Katz = 
\rH^0(X_1(N)'_\FF,\ul{\omega}^{\otimes k}(-\cusps)),
\qquad\text{if $N\geq5$}.
$$
As $\FF$ is flat over~$\FF_p$, we have:
$$
S_k(\Gamma_1(N),\FF)'_\Katz = \FF\otimes S_k(\Gamma_1(N),\FF_p)'_\Katz .
$$
Let $\eps\colon(\ZZ/N\ZZ)^*\to\FF^*$ be a character with $\FF$ a
finite field of characteristic~$p$. For any integer $k$ we have the
$\FF$-vector space $S_k(\Gamma_1(N),\eps,\FF)'_\Katz$ of cuspidal
modular forms on $[\Gamma_1(N)]'$ of weight~$k$, character~$\eps$, and
with coefficients in~$\FF$; it is the subspace of
$S_k(\Gamma_1(N),\FF)'_\Katz$ on which $\ld a\rd$ acts as $\eps(a)$
for all $a$ in~$(\ZZ/N\ZZ)^*$. If $\eps(-1)\neq (-1)^k$ then the space
$S_k(\Gamma_1(N),\eps,\FF)'_\Katz$ is zero.

We have $\FF_p$-linear maps:
$$
S_1(\Gamma_1(N),\FF_p)'_\Katz \lto S_p(\Gamma_1(N),\FF_p)'_\Katz, \quad
A\colon f\mapsto Af,\quad F\colon f\mapsto f^p,
$$
where $A$ is the Hasse invariant. In terms of $q$-expansion at the
standard cusp one has $a_n(Af)=a_n(f)$ and
$a_n(F(f))=a_{n/p}(f)^p=a_{n/p}(f)$ for all $f$ and~$n$, where
$a_{n/p}(f)$ has to be interpreted as zero if $n/p$ is not
integer. (If one works with $[\Gamma_1(N)]$ instead of
$[\Gamma_1(N)]'$ then the $a_n(f)$ lie in $\FF_p(\zeta_N)$ and
$a_{n/p}(f)^p$ is not necessarily equal to~$a_{n/p}(f)$.) It follows
that both $A$ and $F$ are injective.  Combining $A$ and~$F$ gives a
map:
\begin{equation}\label{eqn4.0.1}
\phi\colon {S_1(\Gamma_1(N),\FF_p)'_\Katz}^2\lto
S_p(\Gamma_1(N),\FF_p)'_\Katz, \quad (f,g)\mapsto Af+F(g). 
\end{equation}
A computation using the usual formula
$a_n(T_lf)=a_{nl}(f)+l^{k-1}a_{n/l}(\ld l\rd f)$, valid for $f$ a
modular form of weight~$k>0$, $n$~integer and $l$ any prime (with $\ld
l\rd=0$ if $l$ divides the level), shows that for $a$
in~$(\ZZ/N\ZZ)^*$, and $l\neq p$ prime:
\begin{equation}\label{eq4.0.4}
\ld a\rd\circ\phi = \phi\circ
\begin{pmatrix}
\ld a\rd & 0 \\0 & \ld a\rd
\end{pmatrix},\quad
T_l\circ\phi = \phi\circ
\begin{pmatrix}
T_l & 0 \\0 & T_l
\end{pmatrix},\quad
T_p\circ\phi = \phi\circ
\begin{pmatrix}
T_p & 1 \\-\ld p\rd & 0
\end{pmatrix}.
\end{equation}
The maps obtained from $A$, $F$ and $\phi$ via extension of scalars
via $\FF_p\to \FF$ will be denoted by the same symbols. The
compatibility of $\phi$ with the diamond operators gives an injective
$\FF$-linear map:
\begin{equation}\label{map4.0.5}
F\colon S_1(\Gamma_1(N),\eps,\FF)'_\Katz \lto
S_p(\Gamma_1(N),\eps,\FF)'_\Katz.
\end{equation}
The image of this map $F$ is exactly the subspace of $g$ such that
$a_n(g)=0$ for all $n$ not divisible by~$p$. In order to get an
effective description of this image we will use the derivation
$\theta$ that was constructed by Katz in~\cite{Katz2}. In our
situation, we have an $\FF$-linear map:
$$
\theta\colon S_p(\Gamma_1(N),\eps,\FF)'_\Katz\lto
S_{p+2}(\Gamma_1(N),\eps,\FF)'_\Katz, 
$$
such that for all $g$ in $S_p(\Gamma_1(N),\eps,\FF)'_\Katz$ and all
$n\geq1$ we have:
$$
a_n(\theta(g)) = na_n(g), 
$$
i.e., $\theta$ acts on $q$-expansions as~$qd/dq$. For modular forms of
arbitrary weight, the usual $\theta$ raises the weight by~$p+1$, but
for $g$ of some weight~$pk$, the resulting form of weight $pk+p+1$ has
zeros at the supersingular points, hence can be divided by the Hasse
invariant: see the proof of~\cite[Lemma~3]{Katz2}. The image of $F$
above is the kernel of~$\theta$. The fact that the target of $\theta$
is $S_{p+2}(\Gamma_1(N),\eps,\FF)'_\Katz$ leads to an effective
description of its kernel.

\begin{subprop}\label{prop4.1}
With the notation as above, let $g$ be in
$S_p(\Gamma_1(N),\eps,\FF)'_\Katz$ and put:
$$
B:=\frac{p+2}{12}N\prod_{\substack{l|N\\\text{$l$ prime}}}(1+l^{-1}). 
$$
Suppose that $a_n(g)=0$ for all $n\leq B$ that are not divisible
by~$p$. Then there is a unique $f$ in
$S_1(\Gamma_1(N),\eps,\FF)'_\Katz$ such that $g=F(f)$.
\end{subprop}
\begin{proof}
Let $g$ be as in the statement. Then $\theta(g)$ is a section of
$\ul{\omega}^{\otimes p+2}$ on the stack~$[\Gamma_1(N)]'_\FF$. The
degree of $\omega$ on $[\Gamma_1(1)]'_\FF$ is~$1/24$ (see
\cite[Cor.~10.13.12]{KaMa}, or \cite[VI.4]{DeRa} and use the modular
form $\Delta$ (of weight~12) plus the fact that the automorphism group
of the cusp $\infty$ of $[\Gamma_1(1)]_\FF$ has order~two).  The
degree of $[\Gamma_1(N)]'_\FF$ over $[\Gamma_1(1)]_\FF$ is
$N^2\prod_{l|N}(1-l^{-2})$, where the product is taken over the primes
$l$ dividing~$N$. It follows that the degree of $\ul{\omega}^{\otimes
p+2}$ on the stack~$[\Gamma_1(N)]'_\FF$ is
$(p+2)N^2\prod_{l|N}(1-l^{-2})/24$. It follows that $\theta(g)$, as a
section of~$\ul{\omega}^{\otimes p+2}$, cannot have more zeros than
this degree. We know that $\theta(g)$ has a zero of order at least
$B+1$ at the cusp~$\infty$, and as it is an eigenform for all $\ld
a\rd$ with $a$ in~$(\ZZ/N\ZZ)^*$, the same is true at the cusps in the
$\ZZ/N\ZZ$-orbit of~$\infty$. One checks that these zeros imply that
the divisor of $\theta(g)$ has degree at least~$\phi(N)(B+1)/2$.  (As
we are working on a stack, we have to divide the order of a zero at a
point $x$ by the order of $\Aut(x)$ in order to get the contribution
of $x$ to the degree of the divisor.)  Our choice of $B$ implies that
the divisor of $\theta(g)$ (if non-zero) has degree greater than the
degree of~$\ul{\omega}^{\otimes p+2}$, which means that $\theta(g)=0$.
\end{proof}

\begin{subrem}\label{rem4.2}
In the proof of Proposition~\ref{prop4.1} we have not used that $g$
vanishes at the cusps other than those in the $(\ZZ/N\ZZ)^*$-orbit
of~$\infty$. Hence the result can be improved a little bit, if
necessary. If one knows moreover that $g$ is an eigenform for
Atkin-Lehner (pseudo) involutions, one can replace $B$ by $B/2^r$, where
$r$ is the number of different primes dividing~$N$.  
\end{subrem}
Before we go on, let us record another useful consequence of the
existence of~$\theta$. 
\begin{subprop}\label{prop4.2a}
The map $\phi\colon {S_1(\Gamma_1(N),\FF_p)'_\Katz}^2\lto
S_p(\Gamma_1(N),\FF_p)'_\Katz$ from equation~(\ref{eqn4.0.1}) is
injective.  
\end{subprop}
\begin{proof}
Suppose that $\phi(f,g)=0$. Then we have
$0=\theta(Af+F(g))=\theta(Af)$, hence $0=\theta(f)$. But then $f=0$,
as the weight of $f$ is not divisible by~$p$. It follows that
$F(g)=0$, hence that $g=0$.
\end{proof}
\subsectionrunin{} Now that we know how to characterize the image of
$S_1(\Gamma_1(N),\eps,\FF)'_\Katz$ under~$F$ inside
$S_p(\Gamma_1(N),\eps,\FF)'_\Katz$, it becomes time to investigate how
to compute this ambient vector space. We want to avoid the problems
related to the lifting of elements of
$S_p(\Gamma_1(N),\eps,\FF)'_\Katz$ to characteristic zero with a given
character (see \cite[\S1]{Edix5} and \cite[Prop.~1.10]{Edix5}, they
have to do with what is called Carayol's Lemma). So we describe
$S_p(\Gamma_1(N),\eps,\FF)'_\Katz$ in terms of the $\ZZ[1/N]$-module
$S_p(\Gamma_1(N),\ZZ[1/N])'_\Katz$ of characteristic zero forms of
weight $p$ with no prescribed character.

\begin{subprop}\label{prop4.3}
Let $p$ be a prime number, and $N\geq1$ an integer not divisible
by~$p$. Suppose that $N\neq1$ or that $p\geq5$. Then the map:
$$
\FF_p\otimes S_p(\Gamma_1(N),\ZZ[1/N])'_\Katz \lto 
S_p(\Gamma_1(N),\FF_p)'_\Katz
$$
is an isomorphism, compatible with all Hecke operators $T_n$ (for
$n\geq1$) and $\ld a\rd$ (for $a$ in~$(\ZZ/N\ZZ)^*$). Hence for $\FF$
a finite extension of~$\FF_p$ and a character
$\eps\colon(\ZZ/N\ZZ)^*\to\FF^*$:
$$
S_p(\Gamma_1(N),\eps,\FF)'_\Katz = 
(\FF\otimes S_p(\Gamma_1(N),\ZZ[1/N])'_\Katz)(\eps).
$$
\end{subprop}
\begin{proof}
This follows from \cite[Lemma~1.9]{Edix5} and its proof. 
\end{proof}
\subsectionrunin{}\label{sec4.7}
The next goal is to relate $S_p(\Gamma_1(N),\ZZ[1/N])'_\Katz$ to the
usual description of modular forms as functions on the complex upper
half plane. For $k$ an integer, we let $S_k(\Gamma_1(N),\CC)$ be the
$\CC$-vector space of holomorphic cusp forms on the congruence
subgroup~$\Gamma_1(N)$ of~$\SL_2(\ZZ)$, and, for $A$ a subring
of~$\CC$, we let $S_k(\Gamma_1(N),A)$ be its sub $A$-module of $f$
such that $a_n(f)$ belongs to~$A$ for all~$n$.

With these definitions, $S_k(\Gamma_1(N),\ZZ[1/N])'_\Katz$ and
$S_k(\Gamma_1(N),\ZZ[1/N])$ are the same thing. The following
proposition is well known. It is the final ingredient that we use to
describe the Hecke modules $S_1(\Gamma_1(N),\eps,\FF)'_\Katz$ in terms
of a Hecke algebra that can be computed using group cohomology or
modular symbols.
\begin{subprop}\label{prop4.6}
Let $N\geq1$ and $k$ be integers. Let $\TT$ be the subring of
$\End_\CC(S_k(\Gamma_1(N),\CC))$ generated by the $T_n$,
$n\geq1$. Then $\TT$ contains the $\ld a\rd$ for all $a$
in~$(\ZZ/N\ZZ)^*$, and is generated as $\ZZ$-module by the~$T_n$. The
pairing:
$$
(t,f)\mapsto a_1(tf)
$$
from $\TT\times S_k(\Gamma_1(N),\ZZ)$ to $\ZZ$ is perfect, and gives
an isomorphism of $\TT$-modules:
$$
S_k(\Gamma_1(N),\ZZ) = \TT^\vee.
$$
\end{subprop}
\begin{proof}
For $k\leq0$ the space of cusp forms is zero, and the result is
trivially true. So we suppose that $k\geq1$. The fact that $\TT$
contains all $\ld a\rd$ is well known: one uses that for all prime
numbers $l$ one has $l^{k-1}\ld l\rd = T_l^2-T_{l^2}$, and one takes
two distinct primes $l_1$ and $l_2$ that both have image $a$
in~$\ZZ/N\ZZ$. The fact that the $T_n$ generate $\TT$ as a
$\ZZ$-module follows from the formula that expresses a product
$T_nT_m$ as the sum $\sum_d d^{k-1}\ld d\rd T_{nm/d^2}$ over the
positive common divisors of $n$ and $m$ that are prime to~$N$. To
prove the second statement, one uses the formula $a_n(f)=a_1(T_nf)$
and that $S_k(\Gamma_1(N),\ZZ)$ is by definition the sub $\ZZ$-module
of $f$ in $S_k(\Gamma_1(N),\CC)$ such that $a_n(f)$ is in $\ZZ$ for
all~$n$.
\end{proof}
\begin{subthm}\label{thm4.7}
Let $N\geq5$ be an integer, $p$~a prime number not dividing~$N$,
$\FF$~a finite extension of~$\FF_p$ and
$\eps\colon(\ZZ/N\ZZ)^*\to\FF^*$ a character. Let $\TT$ be the subring
of $\End_\CC(S_p(\Gamma_1(N),\CC))$ generated by the $T_n$,
$n\geq1$. Put $B:=\frac{p+2}{12}N\prod_l(1+l^{-1})$, with the product
taken over the prime divisors of~$N$. Let $L_\eps$ be the sub
$\FF$-vector space of $f$ in $\FF\otimes\TT^\vee=\Hom(\TT,\FF)$ that
satisfy:
\begin{align}
\ld a \rd f& = \eps(a)f &&\text{for all $a$ in $(\ZZ/N\ZZ)^*$;} \\
f(T_n)& = 0  &&\text{for all $n\leq B$ not divisible by~$p$.}
\end{align}
For $f$ in $L_\eps$ one has $f(T_n)=0$ for all $n$ not divisible
by~$p$. For such an~$f$, we have $(T_pf)(T_n)=f(T_{pn})$ for all
$n\geq1$, and $T_p$ induces an isomorphism from $L_\eps$ to its
image~$L'_\eps$. The subspace $L'_\eps$ is stable under all $T_l$ with
$l\neq p$ prime, and by $T_p':=T_p+\eps(p)F$, with $F\colon
L'_\eps\to\FF\otimes\TT^\vee$ defined by: $(Ff)(T_n)=f(T_{n/p})$ (with
$T_{n/p}=0$ if $n/p$ is not integer). Then there is an isomorphism of
Hecke modules between $S_1(\Gamma_1(N),\eps,\FF)'_\Katz$ and $L'_\eps$,
with $T_l$ for $l\neq p$ prime corresponding on both sides, and $T_p$
corresponding to $T_p'$ on~$L'_\eps$.

In terms of $q$-expansions: there is an isomorphism (unique) of
$\FF$-vector spaces between $S_1(\Gamma_1(N),\eps,\FF)'_\Katz$ and
$L'_\eps$ such that for $f$ in $L'_\eps$ its image in
$S_1(\Gamma_1(N),\eps,\FF)'_\Katz$ has $q$-expansion $\sum_n
f(T_n)q^n$.
\end{subthm}
\begin{proof}
Everything follows directly from the previous propositions in this
section, so we just say how those are combined. The image of
$S_1(\Gamma_1(N),\eps,\FF)'_\Katz$ in
$S_1(\Gamma_1(N),\eps,\FF)'_\Katz$ under $F$ is characterized by
Proposition~\ref{prop4.1}. The space
$S_1(\Gamma_1(N),\eps,\FF)'_\Katz$ is described in terms of
$S_p(\Gamma_1(N),\ZZ[1/N])'_\Katz$ by
Proposition~\ref{prop4.3}. In~\S\ref{sec4.7} we have seen that
$S_p(\Gamma_1(N),\ZZ[1/N])'_\Katz$ and $S_p(\Gamma_1(N),\ZZ[1/N])$ are
the same. Proposition~\ref{prop4.6} gives an isomorphism between
$\TT^\vee$ and $S_p(\Gamma_1(N),\ZZ)$.
\end{proof}
It may be of interest to note that twisting modular forms by Dirichlet
characters can be used in order to reduce the number of characters
$\eps$ for which one needs to compute
$S_1(\Gamma_1(N),\eps,\FF)_\Katz$. For example, when $p=2$, all
characters are squares, and twisting can be used to reduce to forms
with trivial character (but possibly a higher level). 

For convenience of the reader we also state a version of
Theorem~\ref{thm4.7} where the character is not fixed. 
\begin{subthm}\label{thm4.8}
Let $N\geq5$ be an integer, and $p$ a prime number not
dividing~$N$. Let $\TT$ be the subring of
$\End_\CC(S_p(\Gamma_1(N),\CC))$ generated by the $T_n$, $n\geq1$. Put
$B':=\frac{p+2}{24}N^2\prod_l(1-l^{-2})$, with the product taken over
the prime divisors of~$N$. Let $L$ be the sub $\FF$-vector space
of $f$ in $\FF\otimes\TT^\vee=\Hom(\TT,\FF)$ that satisfy:
\begin{equation}
f(T_n) = 0 \quad\text{for all $n\leq B'$ not divisible by~$p$.}
\end{equation}
For $f$ in $L$ one has $f(T_n)=0$ for all $n$ not divisible
by~$p$. For such an~$f$, we have $(T_pf)(T_n)=f(T_{pn})$ for all
$n\geq1$, and $T_p$ induces an isomorphism from $L$ to its
image~$L'$. The subspace $L'$ is stable under all $T_l$ with
$l\neq p$ prime, and by $T_p':=T_p+\ld p\rd F$, with $F\colon
L'\to\FF\otimes\TT^\vee$ defined by: $(Ff)(T_n)=f(T_{n/p})$ (with
$T_{n/p}=0$ if $n/p$ is not integer). Then there is an isomorphism of
Hecke modules between $S_1(\Gamma_1(N),\FF_p)_\Katz$ and~$L'$,
with $T_l$ for $l\neq p$ prime corresponding on both sides, and $T_p$
corresponding to $T_p'$ on~$L'$. 

In terms of $q$-expansions: there is an isomorphism (unique) of
$\FF_p$-vector spaces between $S_1(\Gamma_1(N)',\FF_p)_\Katz$
(defined with $\mu_N$-structures) and $L'$ such that for $f$ in
$L'$ its image in $S_1(\Gamma_1(N)',\FF_p)_\Katz$ has
$q$-expansion $\sum_n f'(T_n)q^n$. 
\end{subthm}
\begin{proof}
One just adapts the proof of Theorem~\ref{thm4.7}. The value of $B'$
is gotten by inspection of the proof of Proposition~\ref{prop4.1}.
\end{proof}
For the sake of completeness, we include the following result.
\begin{subprop}\label{prop4.4}
Let $N\geq1$ and $k$ be integers. Then the Hecke operators $T_n$
($n\geq1$) and diamond operators $\ld a\rd$ ($a$ in $(\ZZ/N\ZZ)^*$)
generate the same ring of endomorphisms $\TT$ in
$S_k(\Gamma_1(N),\ZZ[1/N])_\Katz$ and
$S_k(\Gamma_1(N),\ZZ[1/N])'_\Katz$, and these two $\TT$-modules are
isomorphic.
\end{subprop}
\begin{proof}
The description of $[\Gamma_1(N)]'_{\ZZ[1/N]}$ as a twist of
$[\Gamma_1(N)]_{\ZZ[1/N]}$ over $\ZZ[1/N,\zeta_N]$ shows that one of
the two $\TT[1/N]$-modules in question is obtained by twisting the
other over $\ZZ[1/N,\zeta_N]$ via the morphism of groups
$(\ZZ/N\ZZ)^*\to\TT^*$ that sends $a$ to the diamond operator~$\ld
a\rd$. In fact, this twisting can be done on any $\TT[1/N]$-module,
and gives an auto equivalence of the category of
$\TT[1/N]$-modules. Hence it is not surprising that this twisting
functor is isomorphic to the functor $M\otimes_{\TT[1/N]}{-}$, where
$M$ is the twist of $\TT[1/N]$ itself. Just this fact is already
sufficient to conclude that a module and its twist are locally
isomorphic, since $M$ is an invertible $\TT[1/N]$-module. In the
following we will only use that for prime numbers $p$ not
dividing~$N$, $\ZZ_p\otimes S_k(\Gamma_1(N),\ZZ[1/N])_\Katz$ and
$\ZZ_p\otimes S_k(\Gamma_1(N),\ZZ[1/N])$ are isomorphic as
$\ZZ_p\otimes\TT$-modules.

The proof that $M$ is actually free of rank one depends on
computations that we do not want to reproduce here in detail. First
one replaces $\TT[1/N]$ by the group
ring~$\ZZ[1/N][(\ZZ/N\ZZ)^*]$. One has to find a unit $t=\sum_at_a\ld
a\rd$ in $\ZZ[1/N,\zeta_N][(\ZZ/N\ZZ)^*]$ such that:
$$
\sigma_b(t) = \ld b\rd t,\qquad\text{for all $b$ in $(\ZZ/N\ZZ)^*$,}
$$
where $\sigma_b$ is induced by the automorphism of $\ZZ[\zeta_N]$ that
sends $\zeta_N$ to~$\zeta_N^b$. 

Let us just treat the case where $N=p>2$ is prime. Then one can take:
$$
t:=\sum_a(1-\zeta_p^{-a})\ld a^{-1}\rd, \qquad t':=\sum_a\zeta_p^a\ld
a\rd, 
$$
and verify that $tt'=-p$. 

In the case of a prime power one uses the intermediate fields and the
extraction of $p$th roots. The general case is treated by reduction to
the prime power case.
\end{proof}

\begin{subrem}\label{rem4.5}
To complicate things even further, let us mention that the two stacks
$[\Gamma_1(N)]_{\ZZ[1/N]}$ and $[\Gamma_1(N)]'_{\ZZ[1/N]}$ are
isomorphic, via the following construction. To $E/S$ with a closed
immersion $\alpha\colon(\ZZ/N\ZZ)_S\to E$ one associates $E'/S$ and
$\alpha'\colon \mu_{N,S}\to E'$, where $\pi\colon E\to E'$ is the
quotient of $E$ by the image of~$\alpha$, and $\alpha'$ the natural
isomorphism from $\mu_{N,S}$ (the Cartier dual of $(\ZZ/N\ZZ)_S$) to
the kernel of the dual of~$\pi$ (see~\cite[\S2.8]{KaMa}). In the
terminology of~\cite{KaMa}, this isomorphism is called exotic, and
seems to be of no help for comparing $S_k(\Gamma_1(N),\ZZ[1/N])_\Katz$
and $S_k(\Gamma_1(N),\ZZ[1/N])$ as Hecke modules.
\end{subrem}

\section{Some remarks on parabolic cohomology with coefficients modulo~$p$.}
\label{sec4a}
\subsectionrunin{}
Suppose that $N\geq 5$. (Most of what follows still works for
arbitrary $N\geq1$ if one replaces $Y_1(N)$ by the
stack~$[\Gamma_1(N)]_\CC$; we intend to give precise statements in a
future article.)  Let $\pi\colon\EE\to Y_1(N)_\QQ$ be the universal
elliptic curve. Consider the locally constant sheaf
$\rR^1\pi_*\ZZ_\EE$ of free $\ZZ$-modules of rank two
on~$Y_1(N)(\CC)$. For $k\geq2$ let
$\mathcal{F}_k:=\Sym^{k-2}(\rR^1\pi_*\ZZ_\EE)$. Then one defines:
\begin{equation}
\rH^1_\parabolic(Y_1(N)(\CC),\mathcal{F}_k):=
\Image(\rH^1_\compact(Y_1(N)(\CC),\mathcal{F}_k)\to 
\rH^1(Y_1(N)(\CC),\mathcal{F}_k)), 
\end{equation}
with the subscript ``c'' standing for compactly supported
cohomology. Let $j$ denote the inclusion of $Y_1(N)_\QQ$
into~$X_1(N)_\QQ$. Then one has:
\begin{equation}
\begin{split}
\rH^1_\compact(Y_1(N)(\CC),\mathcal{F}_k) & = 
\rH^1(X_1(N)(\CC),j_!\mathcal{F}_k), \\
\rH^1_\parabolic(Y_1(N)(\CC),\mathcal{F}_k) & = 
\rH^1(X_1(N)(\CC),j_*\mathcal{F}_k), \\
\rH^1_\parabolic(Y_1(N)(\CC),\mathcal{F}_2) & = 
\rH^1(X_1(N)(\CC),\ZZ).
\end{split}
\end{equation}
The Shimura isomorphism is the Hodge decomposition
(see~\cite[\S12]{DiIm}):
\begin{equation}
\CC\otimes\rH^1_\parabolic(Y_1(N)(\CC),\mathcal{F}_k) = 
S_k(\Gamma_1(N),\CC)\oplus\ol{S_k(\Gamma_1(N),\CC)}.
\end{equation}
It follows that the Hecke correspondences and diamond operators
generate the same sub $\ZZ$-algebra $\TT(N,k)$ of
$\End(\QQ\otimes\rH^1_\parabolic(Y_1(N)(\CC),\mathcal{F}_k))$ and of
$\End(S_k(\Gamma_1(N),\CC))$. The Hecke modules
$\rH^1_\parabolic(Y_1(N)(\CC),\mathcal{F}_k)$ can be interpreted in
group cohomological terms. To be precise, we mention that:
$$
\rH^1(Y_1(N)(\CC),\mathcal{F}_k) = \rH^1(\Gamma_1(N),\Sym^{k-2}(\ZZ^2)), 
$$
and that the submodule $\rH^1_\parabolic(Y_1(N)(\CC),\mathcal{F}_k)$
is obtained as the subgroup of elements of
$\rH^1(\Gamma_1(N),\Sym^{k-2}(\ZZ^2))$ whose restriction to all
unipotent subgroups of $\Gamma_1(N)$ is zero (see for example
\cite[\S12.2]{DiIm}). For more details on this, especially in relation
to Serre's conjectures in~\cite{Serre1}, one may
consult~\cite{Herremans1}.

For $p$ and $k$ with $p$ a prime not dividing~$N$ and $2\leq k\leq
p+1$, the $\ZZ_p$-module
$\rH^1_\parabolic(Y_1(N)(\CC),\ZZ_p\otimes\mathcal{F}_k)$ is free, and
the morphism:
\begin{equation}\label{eqn4.13.5}
\FF_p\otimes\rH^1_\parabolic(Y_1(N)(\CC),\mathcal{F}_k)\lto
\rH^1_\parabolic(Y_1(N)(\CC),\FF_p\otimes\mathcal{F}_k)
\end{equation}
is an isomorphism. (To prove this, one uses among other things that
the monodromy of $\rR^1\pi_*\ZZ_\EE$ at a cusp is given, up to
conjugation, by some $(\begin{smallmatrix}1 & d\\0 &
1\end{smallmatrix})$, with $d$ dividing~$N$, hence prime to~$p$.) 

The descriptions of spaces of mod~$p$ modular forms of weight one
given in Theorems~\ref{thm4.7} and~\ref{thm4.8} are in terms of the
$\FF_p$-algebra~$\FF_p\otimes\TT(N,p)$. In view of the
isomorphism~(\ref{eqn4.13.5}) one would hope that
$\FF_p\otimes\TT(N,p)$ is the sub $\FF_p$-algebra of
$\End_{\FF_p}(\rH^1_\parabolic(Y_1(N)(\CC),\FF_p\otimes\mathcal{F}_p))$
generated by the $T_n$ and~$\ld a\rd$. But some thinking shows that it
is not clear at all whether or not
$\rH^1_\parabolic(Y_1(N)(\CC),\FF_p\otimes\mathcal{F}_p)$ is a
faithful $\FF_p\otimes\TT(N,p)$-module. As we do know that
$S_k(\Gamma_1(N),\FF_p)_\Katz$ is isomorphic to
$\FF_p\otimes\TT(N,k)^\vee$, when $k\geq2$, some integral $p$-adic
Hodge theory should be applied at this point. The usual theory for
\'etale cohomology and analytic theory give isomorphisms of
$\TT(N,k)$-modules:
\begin{equation}
\rH^1_\parabolic(Y_1(N)(\CC),\FF_p\otimes\mathcal{F}_k) = 
\rH^1_\parabolic(Y_1(N)_{\Qbar_p,\et},\mathcal{F}_{k,p}), 
\end{equation}
where $\mathcal{F}_{k,p}$ is the mod $p$ \'etale counterpart
of~$\mathcal{F}_k$. 
\begin{subthm}\label{thm4.14}
Let $N\geq5$ and $k\geq2$ be integers, and $p$ a prime number not
dividing~$N$. If $k<p$ or $k=p=2$, then
$\rH^1_\parabolic(Y_1(N)(\CC),\FF_p\otimes\mathcal{F}_k)$ is a
faithful $\FF_p\otimes\TT(N,k)$-module. In particular,
$\rH^1(X_1(N)(\CC),\FF_2)$ is a faithful $\FF_2\otimes\TT(N,2)$-module
if $2$ does not divide~$N$. It follows that, for $k<p$ or $k=p=2$, the
$\FF_p$-algebra $\FF_p\otimes\TT(N,k)$ is the sub $\FF_p$-algebra of
endomorphisms of
$\rH^1_\parabolic(Y_1(N)(\CC),\FF_p\otimes\mathcal{F}_k)$ generated by
the $\ld a\rd$ ($a$ in $(\ZZ/N\ZZ)^*$) and $T_n$ ($n\geq1$), and, in
fact, that $\FF_p\otimes\TT(N,k)$ is generated as an $\FF_p$-vector space
by the $\ld a\rd$ ($a$ in $(\ZZ/N\ZZ)^*$) plus the $T_n$, $1\leq n\leq
12^{-1}kN\prod_{l|N}(1+l^{-1})$.
\end{subthm}
\begin{proof}
Let us first deal with the case $2\leq k<p$. We follow Fontaine and
Laffaille as in~\cite{FontaineLaffaille}, adapting it to the special
case that we need, and using some notation as in~\cite{Faltings1}. We
let $\MF_{[0,p-2]}$ be the category of filtered $\phi$-modules
$(D,\Fil,\phi)$ where $D$ is a finite dimensional $\FF_p$-vector
space, $\Fil$~a decreasing filtration on~$D$ such that $D=\Fil^0(D)$
and $\Fil^{p-1}(D)=0$, and $\phi\colon\Gr(D)\to D$ an isomorphism of
$\FF_p$-vector spaces to $D$ from the graded object
$\oplus_{i=0}^{p-2}\Fil^i(D)/\Fil^{i+1}(D)$ associated to
$(D,\Fil)$. Fontaine and Laffaille construct a fully faithful
contravariant functor $\mathbb{V}$ from $\MF_{[0,p-2]}$ to the
category $\Rep_{\FF_p}(G_{\QQ_p})$ of continuous representations of
$G_{\QQ_p}:=\Gal(\Qbar_p/\QQ_p)$ on finite dimensional $\FF_p$-vector
spaces, see~\cite[Thm.~6.1]{FontaineLaffaille}. The objects in the
essential image of $\mathbb{V}$ are called crystalline
representations.

In order to relate this functor $\mathbb{V}$ to our problem, we can
refer to work of Faltings and Jordan (\cite{Faltings1}
and~\cite{FaltingsJordan}), or of Fontaine and Messing and Kato (see
\cite{FontaineMessing}, \cite{Kato1}, or better, the introduction
of~\cite{Breuil1}). In the latter case one has to use the description
of $\rH^1_\parabolic(Y_1(N)(\CC),\FF_p\otimes\mathcal{F}_{k,p})$ as a
suitable piece of $\rH^{k-1}(E^{k-2}(\CC),\FF_p)$, where $E^{k-2}$ is
the proper smooth $\ZZ[1/N]$-model of the $k-2$-fold power of the
universal elliptic curve over~$Y_1(N)$, as described in Deligne
\cite{Deligne1} and Scholl~\cite{Scholl1}. In both cases, one gets the
result that
$\rH^1_\parabolic(Y_1(N)_{\Qbar_p,\et},\mathcal{F}_{k,p})^\vee$ is
crystalline. In the first case, the corresponding object
$(D,\Fil,\phi)$ of $\MF_{[0,p-2]}$ is written
$\rH^1_\parabolic(Y_{\Fbar_p},\Sym^{k-2}(\mathcal{E}))$ (see
\cite[Thm~1.1]{FaltingsJordan}). In the second case, it is a suitable
piece of $\rH^{k-1}_{\deRham}(E^{k-2}_{\FF_p})$.  In both cases,
$(D,\Fil,\phi)$ satisfies (see \cite[Thm~1.1]{FaltingsJordan}):
\begin{equation}
\begin{split}
&\Fil^k(D)= 0\\
&\Fil^{k-1}(D) = \Fil^{1}(D) = S_k(\Gamma_1(N),\FF_p)_\Katz,\\
&\Fil^0(D)/\Fil^1(D) = S_k(\Gamma_1(N),\FF_p)_\Katz^\vee
\end{split}
\end{equation}
By \cite[Thm.~1.2]{FaltingsJordan}, the functor $\mathbb{V}$
transforms the Hecke operators $T_n$ and the $\ld a\rd$ on $D$ into
the duals of the same operators on
$\rH^1_\parabolic(Y_1(N)_{\Qbar_p,\et},\mathcal{F}_{k,p})$. As
$S_k(\Gamma_1(N),\FF_p)_\Katz$ is a faithful
$\FF_p\otimes\TT(N,k)$-module, the same holds for~$D$, hence for
$\rH^1_\parabolic(Y_1(N)_{\Qbar_p,\et},\mathcal{F}_{k,p})$.

Now suppose that $k=p=2$. Let $J$ denote the jacobian of the curve
$X_1(N)$ over~$\ZZ_2$. Then $\rH^1(X_1(N)(\CC),\FF_2)$ is the same as
$\rH^1(X_1(N)(\Qbar),\FF_2)$, and as $\rH^1(X_1(N)(\Qbar_2),\FF_2)$,
after a choice of embeddings of $\Qbar$ into $\CC$
and~$\Qbar_2$. Hence, as $\TT(N,2)$-modules,
$\rH^1(X_1(N)(\CC),\FF_2)$ is the same as $J[2](\Qbar_2)$, the group
of $\Qbar_2$-points of the 2-torsion subgroup scheme of~$J$. Suppose
now that an element $t$ of $\TT(N,2)$ acts as zero on
$\rH^1(X_1(N)(\CC),\FF_2)$. Then it acts as zero on $J[2](\Qbar_2)$,
hence on the group scheme~$J[2]$, since $J[2]_{\QQ_2}$ is reduced and
$J[2]$ is flat over~$\ZZ_2$. Hence $t$ acts as zero on the special
fibre $J[2]_{\FF_2}$ of~$J[2]$. But then $t$ acts as zero on the
tangent space $T_{J[2]_{\FF_2}}(0)$ at~$0$ of~$J[2]_{\FF_2}$. As
$J[2]_{\FF_2}$ contains the kernel of the Frobenius endomorphism
of~$J_{\FF_2}$, we have $T_{J[2]_{\FF_2}}(0)=T_{J_{\FF_2}}(0)$. This
last $\FF_2\otimes\TT(N,2)$-module is known to be free of rank one,
because:
$$
T_{J_{\FF_2}}(0) = \rH^1(X_1(N)_{\FF_2},\calO) =
(S_2(\Gamma_1(N),\FF_2)_\Katz)^\vee.
$$
Hence the image of $t$ in $\FF_2\otimes\TT(N,2)$ is zero. The first
claim of Theorem~\ref{thm4.14} has now been proved.

The claim that $\FF_p\otimes\TT(N,k)$ is the subalgebra of
endomorphisms generated by Hecke correspondences is an immediate
consequence. The last statement then follows from the fact that
$\TT(N,k)$ is generated as $\ZZ$-module by the $\ld a\rd$ ($a$ in
$(\ZZ/N\ZZ)^*$) and the $T_n$ with $1\leq n\leq
12^{-1}kN\prod_{l|N}(1+l^{-1})$, which is proved by an argument
similar to \cite{AgSt2}, see also the proof of
Proposition~\ref{prop4.1}.
\end{proof}
\begin{subrem}
We do not know if
$\rH^1_\parabolic(Y_1(N)(\CC),\FF_p\otimes\mathcal{F}_p)$ is a
faithful $\FF_p\otimes\TT(N,p)$-module, for $N\geq5$ and $p$ a prime
not dividing~$N$. We strongly suspect that this is true and that it
can be proved by the usual reduction to weight two, by working
with~$X_1(Np)$, as for example in \cite{Gross1} and~\cite{Edix4}. Of
course, it would be nicer to have a more direct proof in terms of
$p$-adic Hodge theory. In this case there still is an exact faithful
functor $\mathbb{V}$ from $\MF_{[0,p-1]}$ to $\Rep_{\FF_p}(G_{\QQ_p})$
(see \cite[Thm.~3.3]{FontaineLaffaille}), and
$\rH^{p-1}_{\deRham}(E^{p-2}_{\FF_p})$ is naturally an object of
$\MF_{[0,p-1]}$ (see\cite[Thm.~3.2.1]{BreuilMessing1}). If
$\mathbb{V}$ sends $\rH^{p-1}_{\deRham}(E^{p-2}_{\FF_p})$ to
$\rH^{p-1}(E^{p-2}_{\Qbar_p,\et},\FF_p)$ then it follows that
$\rH^1_\parabolic(Y_1(N)(\CC),\FF_p\otimes\mathcal{F}_p)$ is a
faithful $\FF_p\otimes\TT(N,p)$-module. However it seems not be be
known if $\mathbb{V}$ satisfies this property (see the first sentence
after Theorem~3.2.3 in~\cite{BreuilMessing1}).

We note that for $2\leq k<p$, much more is known about
$\rH^1_\parabolic(Y_1(N)(\CC),\FF_p\otimes\mathcal{F}_k)$ as
$\FF_p\otimes\TT(N,k)$-module, and about the ring $\TT(N,k)$ itself,
see for example \cite[Thm.~2.1]{FaltingsJordan}. Especially, after
localisation at a maximal ideal $m$ of $\FF_p\otimes\TT(N,k)$ such
that the corresponding Galois representation to
$\GL_2((\FF_p\otimes\TT(N,k))/m)$ is irreducible, one knows that
$(\FF_p\otimes\TT(N,k))_m$ is Gorenstein, and with
$\rH^1_\parabolic(Y_1(N)(\CC),\FF_p\otimes\mathcal{F}_k)_m$ free of
rank two over it. In a lot of cases one even knows that
$(\FF_p\otimes\TT(N,k))_m$ is a complete intersection, see
\cite{Diamond1}, and also~\cite{Dickinson1}.
\end{subrem}

\subsectionrunin{} The most interesting case for the study weight one
forms mod~$p$ is when $p=2$. Firstly, the exceptional case in Serre's
conjectures has not been proved and should be subject to some
experimental checking. To be precise, let $\rho\colon
G_\QQ\to\GL_2(\Fbar_2)$ be a continuous irreducible representation
that is modular of some type, and unramified at~$2$. If
$\rho(\Frob_2)$ has a double eigenvalue, then it is not known whether
$\rho$ comes from a form of weight one (see the introduction
of~\cite{Edix4}). If $\rho(\Frob_2)$ is scalar, then it is not known
whether $\rho$ comes from the expected level (see \cite{Buzzard1} and
Theorem~3.2 in~\cite{Buzzard2}). Secondly, when $p=2$ there is no
distinction between even and odd Galois representations mod~$p$, hence
mod~$2$ modular forms conjecturally lead to \emph{all} continuous
irreducible representations from $G_\QQ$ to~$\GL_2(\Fbar_2)$. In this
case, one may want to investigate the spaces
$S_1(\Gamma_0(N),\FF_2)_\Katz$ of modular forms on~$\Gamma_0(N)$, with
$N$~odd. (Note that the character of a non-zero mod~$p$ form of weight
one is necessarily non-trivial if $p\neq2$.)  Working directly with
$\Gamma_0(N)$ has the advantage of being computationally faster, and
maybe easier. So it is useful to know the relations between the spaces
$S_2(\Gamma_0(N),\FF_2)_\Katz$ (that give us the
$S_1(\Gamma_0(N),\FF_2)_\Katz$) and $\rH^1(\Gamma_0(N),\FF_2)$ (that
one can compute). As usual, such relations are given via their
algebras of endomorphisms generated by Hecke operators.
\begin{subthm}
Let $N\geq 5$ be odd. Then $S_2(\Gamma_0(N),\FF_2)_\Katz$ and
$\rH^1_\parabolic(\Gamma_0(N),\FF_2)$ have the same non-Eisenstein
systems of eigenvalues for all~$T_n$, $n\geq1$. (A system of
eigenvalues is called Eisenstein if its associated Galois
representation is reducible.)
\end{subthm}
\begin{proof}
Let $\TT$ denote $\TT(N,2)$, i.e., the ring of endomorphisms of
$S_2(\Gamma_1(N),\CC)$ generated by the~$T_n$ ($n\geq1$) and the~$\ld
a\rd$ ($a$~in~$(\ZZ/N\ZZ)^*$). Let $I$ be the ideal in $\TT$ generated
by the $\ld a\rd -1$ ($a$~in~$(\ZZ/N\ZZ)^*$). Then, as $N\geq5$, we
have $S_2(\Gamma_1(N),\FF_2)_\Katz = \FF_2\otimes\TT^\vee$, and:
$$
S_2(\Gamma_0(N),\FF_2)_\Katz =
S_2(\Gamma_1(N),\FF_2)_\Katz^{(\ZZ/N\ZZ)^*} =
(\FF_2\otimes(\TT/I))^\vee.
$$
It follows that the algebra of endomorphisms of
$S_2(\Gamma_0(N),\FF_2)_\Katz$ generated by the~$T_n$ ($n\geq1$) is
$\FF_2\otimes(\TT/I)$. 

On the group cohomology side we have an exact sequence (coming from
the Hochschild-Serre spectral sequence), with
$\Delta:=(\ZZ/N\ZZ)^*$:
$$
\rH^1(\Delta,\FF_2) \to \rH^1(\Gamma_0(N),\FF_2) \to
\rH^1(\Gamma_1(N),\FF_2)^\Delta \to \rH^2(\Delta,\FF_2).
$$
As the outer two terms are Eisenstein, we get an isomorphism between
the middle two after localising away from the Eisenstein maximal
ideals. Hence the natural map from
$\rH^1_\parabolic(\Gamma_0(N),\FF_2)$ to
$\rH^1(\Gamma_1(N),\FF_2)^\Delta$ becomes an isomorphism, after
localising away from the Eisenstein maximal
ideals. Theorem~\ref{thm4.14} says that the algebra of endomorphisms
of $\rH^1_\parabolic(\Gamma_1(N),\FF_2)$ is a faithful
$\FF_2\otimes\TT$-module. It follows that the support of
$\rH^1_\parabolic(\Gamma_1(N),\FF_2)^\Delta$ in
$\Spec(\FF_2\otimes\TT)$ (as a set of prime ideals) is that of
$(\TT/I)\otimes_\TT\rH^1_\parabolic(\Gamma_1(N),\FF_2)^\vee$, hence
equal to $\Spec(\FF_2\otimes(\TT/I))$. This proves that
$S_2(\Gamma_0(N),\FF_2)_\Katz$ and
$\rH^1_\parabolic(\Gamma_0(N),\FF_2)$ have the same non-Eisenstein
systems of eigenvalues.
\end{proof}
\begin{subthm}
Let $N\geq5$ be odd and divisible by a prime number $q\equiv-1$
modulo~$4$ (hence the stabilizers of the group $\Gamma_0(N)/\{1,-1\}$
acting on~$\HH$ have odd order). Then $S_2(\Gamma_0(N),\FF_2)_\Katz$
and $\FF_2\otimes S_2(\Gamma_0(N),\ZZ)$ are equal, and the
localisations at non Eisenstein maximal ideals of the algebras of
endomorphisms of $S_2(\Gamma_0(N),\FF_2)_\Katz$ and
$\rH^1_\parabolic(\Gamma_0(N),\FF_2)$ generated by all $T_n$
($n\geq1$) coincide: both are equal to that of $S_2(\Gamma_0(N),\ZZ)$
tensored with~$\FF_2$.
\end{subthm}
\begin{proof}
We keep the notations of the previous proof. We start working on the
modular forms side. As $X_1(N)_{\FF_2}\to X_0(N)_{\FF_2}$ is not
wildly ramified, we have:
$$
\rH^0(X_0(N)_{\FF_2},\Omega) = \rH^0(X_1(N)_{\FF_2},\Omega)^\Delta = 
S_2(\Gamma_1(N),\FF_2)_\Katz^\Delta = S_2(\Gamma_0(N),\FF_2)_\Katz .
$$
In the previous proof we have seen that $S_2(\Gamma_0(N),\FF_2)_\Katz
= (\FF_2\otimes(\TT/I))^\vee$, hence we see that the dimension
over~$\FF_2$ of $\FF_2\otimes(\TT/I)$ is the genus
of~$X_0(N)_{\FF_2}$. On the other hand we have:
$$
S_2(\Gamma_0(N),\CC) = S_2(\Gamma_1(N),\CC)^\Delta =
\CC\otimes(\TT/I)^\vee. 
$$
It follows that $\ZZ_{(2)}\otimes\TT/I$ is a free
$\ZZ_{(2)}$-module. Let $\TT_0$ be the image of $\TT$ in the ring of
endomorphisms of $S_2(\Gamma_0(N),\ZZ)$. As
$S_2(\Gamma_0(N),\ZZ)=(\TT/I)^\vee$, we see that $\TT_0$ is the
quotient of $\TT/I$ by its ideal of torsion elements. It follows that
the natural map $\TT/I\to\TT_0$ becomes an isomorphism after tensoring
with~$\ZZ_{(2)}$. We have now shown that the Hecke algebra of
$S_2(\Gamma_0(N),\FF_2)_\Katz$ is~$\FF_2\otimes\TT_0$. 

Let us now consider the group cohomology side. The action of
$\Gamma_0(N)$ on~$\HH$ factors through the quotient
$\Gamma_0(N)/\{1,-1\}$, and we have an exact sequence:
$$
0\to \rH^1(\Gamma_0(N)/\{1,-1\},\FF_2) \to \rH^1(\Gamma_0(N),\FF_2)
\to \rH^1(\{1,-1\},\FF_2) = \FF_2,
$$
with the last term a module with Eisenstein Hecke eigenvalues.  The
stabilizers for the action of $\Gamma_0(N)/\{1,-1\}$ on~$\HH$ are of
odd order, hence the natural map from $\rH^1(Y_0(N)(\CC),\FF_2)$ to
$\rH^1(\Gamma_0(N)/\{1,-1\},\FF_2)$ is an isomorphism.  Hence, after
localising away from the Eisenstein maximal ideals,
$\rH^1_\parabolic(\Gamma_0(N),\FF_2)$ and $\rH^1(X_0(N)(\CC),\FF_2)$
are the same. But the latter is a faithful $\FF_2\otimes\TT_0$-module,
by the arguments in the proof of Theorem~\ref{thm4.14}.
\end{proof}

\section{Eigenspaces in weight $p$ and eigenforms of weight one.}
\label{sec5}
\subsectionrunin{} 
The results of Section~\ref{sec4} make it possible to compute the
spaces $S_1(\Gamma_1(N),\FF_p)_\Katz$, with $N\geq5$ and $p$ a prime
not dividing~$N$, as Hecke modules. This is useful if one is
interested in ring theoretic properties of the associated Hecke
algebras. Sometimes, one is just interested in the systems of
eigenvalues. The aim of this section is to give a result about
eigenspaces.

\begin{subprop}\label{prop5.1}
Let $N\geq1$ be an integer, $p$ a prime number that does not
divide~$N$, and $\eps\colon (\ZZ/N\ZZ)^*\to\Fbar_p^*$ a character. Let
$V$ be a common non-zero eigenspace in
$S_p(\Gamma_1(N),\eps,\Fbar_p)_\Katz$ for all~$T_l$, $l\neq p$
prime. Let $a$ denote the system of eigenvalues attached to~$V$: $T_l$
acts as~$a_l$ on~$V$. Then $V$ is of dimension at most two. If
$\dim(V)=1$, then $a$ does not occur in
$S_1(\Gamma_1(N),\eps,\Fbar_p)_\Katz$. If $\dim(V)=2$, then the
eigenspace $V_1$ in $S_1(\Gamma_1(N),\eps,\Fbar_p)_\Katz$ attached to
$a$ is of dimension one. In this case, the $T_p$ eigenvalue of $V_1$
is the trace of $T_p$ on~$V$.
\end{subprop}
\begin{proof}
We consider $q$-expansions at the standard cusp, and use the formulas
for the action of the $T_n$ and~$\ld a\rd$. Let $V$ be as in the
proposition. Let $d$ be the dimension of~$V$. As $T_p$ commutes with
all~$T_l$, $T_p$ acts on~$V$, hence $V$ contains a non-zero
eigenvector $f$ for all~$T_n$ ($n\geq1$). In particular,
$a_1(f)\neq0$. It follows that the subspace $V'$ of $V$ consisting of
the $g$ in $V$ with $a_1(g)=0$ is of dimension~$d-1$. For every
element $g$ of~$V'$, we have $0=a_1(T_n(g))=a_n(g)$ for all $n$ not
divisible by~$p$. It follows that $V'=FV_1$, with $F$ as
in~(\ref{map4.0.5}), and $V_1$ the eigenspace of
$S_1(\Gamma_1(N),\eps,\Fbar_p)_\Katz$ associated to~$a$. 

Suppose that $g$ is in $V_1$ and that $a_1(g)=0$. Then $a_n(g)=0$ for
all $n$ prime to~$p$. But then $\theta(g)=0$, implying that $g=0$, as
$p$ does not divide the weight of~$g$ (see~\cite{Katz2}). We have
proved that $d\leq2$. 

If $V_1\neq\{0\}$, then $d\geq 2$ by Proposition~\ref{prop4.2a} and
the identities~(\ref{eq4.0.4}). Hence if $d=1$, then $V_1=\{0\}$, and
$a$ does not occur in $S_1(\Gamma_1(N),\eps,\Fbar_p)_\Katz$. If $d=2$,
then $\dim(V_1)=1$, and the trace of $T_p$ on $V$ is the
$T_p$-eigenvalue $a_p$ of $V_1$ by the identities~(\ref{eq4.0.4}).
\end{proof}
\begin{subrem}
For computational purposes, let us add the following. In
Proposition~\ref{prop5.1} one only needs to take into account the
$T_l$ with $l\leq B$, where $B$ is as in
Proposition~\ref{prop4.1}. The proof that these $T_l$ define the same
eigenspaces as all the $T_l$ with $l\neq p$ is given by
Proposition~\ref{prop4.1}.
\end{subrem}

\newpage
\appendix
\section{Lettre de Mestre \`a Serre}
\begin{center}Par Jean-Fran{\c c}ois Mestre\footnote{
Address: 
Centre de Math{\'e}matiques de Jussieu,
Projet Th{\'e}orie des Nombres,
Universit{\'e} Paris~7,
Etage~9, bureau~9E10,
175, rue de Chevaleret,
75013 PARIS, 
France;
E-mail: {\tt mestre@math.jussieu.fr}
}
\end{center}

\begin{flushright}
Paris, 8 Octobre 1987
\end{flushright}
\vspace{10ex}

\begin{center}
Cher Monsieur,
\end{center}
\vspace{10ex}

Voici quelques r\'esultats concernant les formes $\modulo 2$ de
poids~$1$, en niveau premier.

On utilise votre m\'ethode pour trouver une base des formes $\modulo
2$, niveau~$p$, poids~$1$, qui, {\'e}lev{\'e}es {\`a} la
puissance~$2$, sont la r\'eduction de formes sur~$\CC$ de niveau~$p$,
poids~$2$.

Tous les conducteurs $p$ premiers $\leq 1429$ ont \'et\'e \'etudi\'es.

\subsection{Repr\'esentations di\'edrales, de type~$S_4$ et de 
type~$A_5$}

Soit $k$ l'extension quadratique non ramifi\'ee en dehors de~$p$, de
nombre de classes $h=2^mu$, avec $u$ impair.

Vous sugg\'erez qu'on devrait alors trouver $(u-1)/2$ formes
modulaires $\modulo 2$, de poids~$1$ et conducteur~$p$, correspondant
\`a la repr\'esentation di\'edrale associ\'ee (que $k$ soit r\'eel ou
imaginaire).

Dans tous les cas examin\'es, (i.e. $p\leq 1429$), on a toujours
trouv\'e les formes en question (plus pr\'ecis\'ement, $(u-1)/2$
formes de poids~$1$ niveau~$p$, \`a coefficients dans~$\FF_q$,
$q=2^{(u-1)/2}$, conjugu\'ees entre elles, dont les premiers
coefficients $a_l$ sont compatibles avec ceux de la repr\'esentation
di\'edrale).

\medskip
Parfois, l'espace des formes de poids~$1$ est de dimension
strictement plus grande que le nombre de syst\`emes de valeurs propres
des~$T_l$.

Disons que, pour un conducteur~$p$, on est dans le cas~$B(m)$ si un
op\'erateur~$T_l$ a une valeur propre~$a_l$ telle que $T_l-a_l$ est
nilpotent de degr\'e $m\geq 2$.

Le tableau suivant indique les valeurs de~$p$ ``exceptionnelles''
(i.e. pour lesquelles on est dans le cas $B(m)$ ($m\geq 2$), ou bien
o\`u on a trouv\'e des formes ne correspondant pas \`a des
repr\'esentations di\'edrales).

$$\begin{array}{ccccccccccc}
p&229&257&283&331&491&563&643&653&751\\
d&2&2&3&3&6&6&3&4&9\\
h&3&3&3&3&9&9&3&1&15\\
 & & & & & & & & & \\
p&761&1061&1129&1229&1367&1381&1399&1423&1429\\
d&2&4&5&2&16&4&15&6&8\\
h&3&1&9&3&25&1&27&9&5
\end{array}$$
o\`u $h$ est le nombre de classes du corps quadratique (r\'eel ou
imaginaire suivant que $p\equiv 1$ ou $3\modulo 4$) non ramifi\'e en
dehors de $p$, et o\`u $d$ est la dimension de l'espace des formes de
poids $1\modulo 2$ trouv\'ees.

\begin{itemize}
\item
Pour $p=229$, $283$,
$331$, $491$, $563$, $643$, $751$, $1399$ et $1423$,
on est dans le cas $B(m)$,
avec $m=2$ pour $229$ et $m=3$ sinon.
Dans ces divers cas, les tables de Godwin montrent l'existence
d'une extension de type $S_4$ non ramifi\'ee en dehors de $p$,
ce qui pourrait expliquer le ph\'enom\`ene d'unipotence constat\'e.

\item
Pour $p=257$, $761$, $1129$ et $1229$, on est dans le cas $B(2)$.

\item
Pour $p=1367$, on est dans le cas~$B(3)$: le groupe des classes
d'id\'eaux de $k=\QQ(\sqrt{-p})$ est isomorphe \`a~$\ZZ/25\ZZ$. \`A la
repr\'esentation di\'edrale de degr\'e~$10$ de $\QQ$ associ\'ee
correspondent $2$ syst\`emes conjugu\'es de valeurs propres dans
$\FF_4$.

L'espace primaire associ\'e \`a l'un quelconque de ces $2$ syst\`emes
est de dimension~$3$.
\item
Le cas $653$ (resp. $1381$) correspond sans doute (d'apr\`es l'examen
des $a_l$ pour $l$ petit) au corps de type $A_5$ engendr\'e par les
racines de $x^5+3x^3+6x^2+2x+1=0$
(resp. $x^5+3x^4+10x^3+4x^2-16x-48=0$). Ces deux corps sont non
ramifi\'es en $2$, ce qui est en accord avec votre conjecture.
\end{itemize}

Dans tous les cas ci-dessus, il existe sur $\CC$ une forme de poids $1$,
mais de niveau \'eventuellement multiple strict de $p$, qui $\modulo 2$
donne la forme obtenue.

\subsection{Formes de poids $1\modulo 2$ ne provenant pas de formes de poids
$1$ en caract\'eristique $0$}

Pour $p=1429$, la dimension des formes de poids $1\modulo 2$ est
de dimension $8$.

Les polyn\^omes caract\'eristiques des op\'erateurs de Hecke $T_l$
sont, pour $l$ premier $\leq 13$:

$$\begin{array}{ccc}
2&3&5\\
x^2(x^3+x^2+1)^2&(x^2+x+1)(x^3+x^2+1)^2&(x^2+x+1)(x^3+x+1)^2\\
&&\\
7&11&13\\
(x^2+x+1)(x^3+x+1)^2&x^2(x^3+x+1)^2&(x^2+x+1)(x^3+x^2+1)^2
\end{array}$$

Le nombre de classes du corps quadratique r\'eel $\QQ(\sqrt{1429})$,
(non ramifi\'e en dehors de $1429$) est $5$. Ceci explique (modulo
votre conjecture) les $2$
formes \`a coefficients dans $\FF_4$, propres pour les op\'erateurs de
Hecke, dont l'existence est mise en \'evidence par les \'equations
ci-dessus.

Il reste trois formes propres, \`a coefficients conjugu\'es dans $\FF_8$.
Si $r\in \FF_8$ v\'erifie l'\'equation $r^3+r^2+1=0$, l'une d'elle a comme
coefficients $a_p$, pour $p$ premier $\leq 89$:

$$\begin{array}{cccccccccc}
2&3&5&7&11&13&17&19&23&29\\
r&r&r^2+r&r^2+1&r+1&r&r^2&0&r^2+r+1&r^2\\
& & & & & & & & &\\
31&37&41&43&47&53&59&61&67&71\\
r^2+r+1&r^2+r&r^2+r+1&r^2&r&r^2+r&r^2+1&1&r^2+r&r+1\\
& & & & & & & & &\\
73&79&83&89&&&&&&\\
r&r^2+r&1&1&&&&&&
\end{array}$$

Le fait que l'on obtient $1$ et $r$ comme traces de Frobenius
montre que l'image de $G_Q$ dans $\SL_2(\FF_8)$ (pour la
repr\'esentation correspondant \`a la forme en question)
est $\SL_2(\FF_8)$ tout entier.

{\bf Par suite, on ne peut obtenir cette
repr\'esentation \`a partir d'une forme de poids $1$ sur $\CC$.}

En poussant plus loin les calculs, on trouve quatre
autres cas similaires, donnant des re\-pr\'e\-sen\-ta\-tions de type
$\SL_2(\FF_8)$ (obtenus pour les conducteurs $p=1613$, $1693$, $2017$ et
$2089$).

\vspace{10ex}

Bien \`a vous,

\medskip
\begin{center}
J-F. Mestre.
\end{center}

\section{Computing Hecke algebras of weight $1$ in MAGMA}
\addtocontents{toc}{(by Gabor Wiese)}
\begin{center}By Gabor Wiese\footnote{Supported by the European
Research Training Network Contract HPRN-CT-2000-00120 ``Arithmetic
Algebraic Geometry''. 
Address: 
Mathematisch Instituut,
Universiteit Leiden,
Postbus 9512,
2300 RA Leiden,
The Netherlands;
{\tt http://www.math.leidenuniv.nl/$\sim$gabor/}, e-mail:
{\tt gabor@math.leidenuniv.nl}}
\end{center}

\subsection{Introduction}

The aim of this appendix is twofold.  On the one hand, we report on an
implementation in MAGMA (see \cite{Magma}) of a module for the Hecke
algebra of Katz cusp forms of weight $1$ over finite fields, which is
based on section 4 of this article.

On the other hand, we present results of computations done in relation
with the calculations performed by Mestre (see appendix A) in 1987.

The program consists of two packages, called {\em Hecke1} and {\em
CommMatAlg}.  The source files and accompanying documentation
(\cite{hdoc} and \cite{cma}) can be downloaded from the author's
homepage ({\tt http://www.math.leidenuniv.nl/$\sim$gabor/}).

The author would like to express his gratitude to Bas Edixhoven for
his constant support.

\subsection{Algorithm}

In the current release MAGMA (\cite{Magma}) provides 
William Stein's package HECKE, which contains 
functions for the computation of Hecke algebras and modular forms over
fields. There is, however, the conceptual restriction to
{\em weights greater equal $2$}.

Edixhoven's approach for the construction of a good weight $1$ Hecke
module, which is at the base of the implemented algorithm, relates the
Hecke algebra of characteristic $p$ Katz cusp forms of weight $1$ to
the Hecke algebra of classical weight $p$ cusp forms over the complex
numbers.  The latter can for instance be obtained using modular
symbols.

\subsubsection*{Katz modular forms}

Following the notations of section 4, we denote by $\ksfc k 1 N
{\overline{\epsilon}} \FF$ the space of {\em Katz cusp forms} of
weight $k$, level $N$, with character $\overline{\epsilon}: (\ZZ/N)^*
\to \FF^*$ over the $\ZZ[1/N]$-algebra $R$, where we impose that $k
\ge 1$ and $N \ge 5$.  For a definition see section 4 or \cite{Gross1}
for more details.

By the space of {\em classical cusp forms} $\Sfc k 1 N \epsilon R$
over a ring $R \subseteq \CC$, we understand the sub-$R$-module of
$\Sfc k 1 N \epsilon \CC$ consisting of the forms with Fourier
coefficients (at infinity) in the ring $R$.

Let us mention that for a homomorphism of $\ZZ[1/N]$-algebras $R \to
S$, we have the isomorphism (\cite{Gross1}, Prop.~2.5, and the proof
of \cite{DiIm}, Thm.~12.3.2)
$$ \ksf k 1 N R \otimes_R S \cong \ksf k 1 N S,$$
if $k \ge 2$ or if $R \to S$ is flat. 
Using the statements in 4.7, it follows in particular that we have the
equality 
$$ \ksf k 1 N R = \Sf k 1 N R,$$
in case that $\ZZ[1/N] \subseteq R \subseteq \CC$ or $k \ge 2$.

\subsubsection*{Modular symbols}

Given integers $k \ge 2$ and $N \ge 1$, one can define the complex
vector space $\ms k 1 N$ of {\em cuspidal modular symbols} (see
e.g.~\cite{Merel1}, section 1.4).  On it one has in a natural manner
Hecke and diamond operators, and there is a non-degenerate pairing
$$ \big( \Sf k 1 N \CC \oplus \overline{\Sf k 1 N \CC } \big) \times
\ms k 1 N \to \CC,$$
with respect to which the diamond and Hecke operators are adjoint
(see \cite{Merel1}, Thm.~3 and Prop.~10).

We recall that the diamond operators provide a group action of
$(\ZZ/N)^*$ on the above spaces.  For a character $\epsilon: (\ZZ/N)^*
\to \CC^*$ one lets, in analogy to the modular forms case, $\msc k 1 n
\epsilon$ be the $\epsilon$-eigenspace.

Let $\ZZ[\epsilon]$ be the smallest subring of $\CC$ containing all
values of $\epsilon$. It follows that the $\ZZ[\epsilon]$-algebra
generated by all Hecke operators acting on $\Sfc k 1 N \epsilon \CC$
is isomorphic to the one generated by the Hecke action on $\msc k 1 n
\epsilon$.  The same applies to the $\ZZ$-algebra generated by the
Hecke operators on the full spaces (i.e.~without a character).

\begin{subsubnot}\label{not}
We call the Hecke algebras described here above $\TT(\epsilon)$ and
$\TT$ respectively.
\end{subsubnot}

It is known (for the method see e.g.~Prop.~4.2) that the first $Bk$
Hecke operators suffice to generate $\TT(\epsilon)$, where the number
$B$ is $\frac{N}{12} \prod_{l \mid N, l \text{ prime}} (1 +
\frac{1}{l})$.  For the full Hecke algebra $\TT$ one has to take $B k
\varphi(N) / 2$.

\subsubsection*{Weight $1$ as subspace in weight $p$}

Let us assume the following
\begin{subsubsetting}\label{setting}
Let $K$ be a number field, $\cO_K$ its ring of integers, $\fP$ a prime of
$\cO_K$ above the rational prime $p$ and $N\ge 5$ an integer coprime
to~$p$. Moreover, we consider a character 
$\epsilon: (\ZZ/N)^* \to \cO_K^*$.
For a given field extension $\FF$ of $\cO_K / \fP$, we fix the canonical 
ring homomorphism
$\phi: \cO_K \twoheadrightarrow \cO_K / \fP \hookrightarrow \FF$.
We denote by $\overline{\epsilon}$ the composition of $\epsilon$ with
$\phi$.
Recall that $B$ was defined to be 
$\frac{N}{12} \prod_{l \mid N, l \text{ prime}} (1 + \frac{1}{l})$.
\end{subsubsetting}

We shall quickly explain how Edixhoven relates weight $1$ to weight $p$
in section 4 in order to be able to formulate our statements.

The main tool is the {\em Frobenius} homomorphism $F: \ksf 1 1 N
{\FF_p} \to \ksf p 1 N {\FF_p}$ defined by raising to the $p$-th
power.  Hence on $q$-expansions it acts as $a_n(Ff) = a_{n/p}(f)$,
where $a_{n/p}(f) = 0$ if $p \nmid n$.  Also by $F$ we shall denote
the homomorphism obtained by base extension to $\FF$.
One checks that $F$ is compatible with the character. 
The sequence of $\FF$-vector spaces
$$ 0 \to  \ksfc 1 1 N {\overline{\epsilon}} \FF 
     \xrightarrow{F} \ksfc p 1 N {\overline{\epsilon}} \FF 
     \xrightarrow{\Theta} \ksfc {p+2} 1 N {\overline{\epsilon}} \FF$$
is exact, where $\Theta$ denotes the derivation described before
Prop.~4.2.  The image of $F$ in $\ksfc p 1 N {\overline{\epsilon}} \FF
$ is effectively described by Prop.~4.2 to be those $f \in \ksfc p 1 N
{\overline{\epsilon}} \FF$ such that $a_n(f) = 0$ for all $n$ with $p
\nmid n$, where it suffices to take $n \le B(p+2)$ with $B$ as before.

Using the homorphisms
\begin{align*}
\ksfc 1 1 N {\overline{\epsilon}} \FF & 
  \overset{F}{\hookrightarrow} \ksfc p 1 N {\overline{\epsilon}} \FF 
  \overset{\text{4.6}}{\cong} 
   \big( (\Sf p 1 N \ZZ) \otimes_\ZZ \FF \big) (\overline{\epsilon})\\
& \overset{\text{4.8}}{\cong} 
   \big( \Hom_\ZZ (\TT,\ZZ) \otimes_\ZZ \FF) \big) (\overline{\epsilon})
  \cong (\TT \otimes_\ZZ \FF)^\vee (\overline{\epsilon}),
\end{align*}
one obtains an isomorphism of Hecke modules (cp.~Thm.~4.9)
\begin{equation} \label{eqthm1}
        \ksfc 1 1 N {\overline{\epsilon}} \FF 
        \cong \big( (\TT \otimes_\ZZ \FF) / \tilde{\cR} \big)^\vee, 
\end{equation}
where $\tilde{\cR}$ denotes the sub-$\FF$-vector space of $\TT
\otimes_\ZZ \FF$ 
generated by $1 \otimes \overline{\epsilon}(l) - <l> \otimes 1$ 
for $(l,N)=1$ and by $T_n$ for $n \le B(p+2)$ and $p \nmid n$.
The action of the Hecke operators is the same as the one given in the
proposition below.

We would like to replace the full Hecke algebra $\TT$, which is
expensive to calculate,  
by $\TT(\epsilon)$. One has a natural surjection
$ \TT \otimes_\ZZ \ZZ[\epsilon] \twoheadrightarrow \TT(\epsilon)$, which
sends $<l> \otimes 1$ to $\epsilon(l) \cdot \id$. 

\begin{subsubprop}\label{thm2}
Assume the setting \ref{setting} and the notation \ref{not}.
Let $\cR$ be the sub-$\FF$-vector space of 
$\TT(\epsilon) \otimes_{\ZZ[\epsilon]} \FF$ generated by 
$T_n \otimes 1$ for those $n \le B(p+2)$ not divisible by $p$.
Then there is an injection of Hecke modules
$$ \big( (\TT(\epsilon) \otimes_{\ZZ[\epsilon]} \FF) / \cR \big)^\vee
   \hookrightarrow   \ksfc 1 1 N {\overline{\epsilon}} \FF.$$

For a prime $l \neq p$, the natural action of the Hecke operator $T_l$
in weight $p$ corresponds to the action of $T_l$ in weight $1$.  The
natural action of the operator $T_p + \overline{\epsilon}(p) F$ on the
left corresponds to the action of $T_p$ in weight $1$.  Here
$F: \TT(\epsilon) \otimes_{\ZZ[\epsilon]} \FF \to 
    \TT(\epsilon) \otimes_{\ZZ[\epsilon]} \FF$ 
sends $T_n \otimes 1$ to $T_{n/p} \otimes 1$ with the
convention $T_{n/p} \otimes 1 =0$ if $p$ does not divide $n$.
\end{subsubprop}

\begin{proof}
With $\TT$ and $\tilde{\cR}$ as defined before the proposition, we
have a surjection
$$ (\TT \otimes_\ZZ \FF)/\tilde{\cR} \twoheadrightarrow 
   (\TT(\epsilon) \otimes_{\ZZ[\epsilon]} \FF) / \cR.$$
Now taking $\FF$-vector space duals together with equation
\ref{eqthm1} gives the claimed injection. The explicit form of the
operators follows immediately from equation 4.1.2.
\end{proof}

We treat a special case separately.
\begin{subsubcor}
Take in propostion \ref{thm2} the trivial character $1$ and $p = 2$. Then
the injection is an isomorphism if there is a prime $q$ dividing $N$ such
that $q \equiv 3$ modulo $4$.
\end{subsubcor}

\begin{proof}
As in the proof of Thm.~5.6, one shows that the Hecke algebra
of $\ksf 2 0 N {\FF_2}$ is $\TT(1) \otimes \FF_2$. Hence, we have
$$ \big(\TT(1) \otimes \FF_2\big)^\vee \cong 
   \big(\TT \otimes \FF_2 / 
(1\otimes 1 - <l> \otimes 1 \; | \; (l,N)=1 ) \big)^\vee,$$
whence the corollary follows.
\end{proof}

\subsection{Software}

\subsubsection*{Functionality}

In this section we wish to present, in a special case, what {\em
Hecke1} computes. Please consult section $2$ of \cite{hdoc} for
precise statements.

\noindent{\bf \underline{INPUT:}} 
Let $C$ be the space of cuspidal modular symbols of weight $p = 2$ and
odd level $N \ge 5 $ for the trivial character over the rational
numbers.

\noindent{\bf \underline{COMPUTE:}} 
Let $\phi: \ZZ \to \FF_2$ be the canonical ring homomorphism. We
denote by $\overline{T_i}$ the image under $\phi$ of the matrix
representing the $i$-th Hecke operator $T_i$ acting on the natural
integral structure of $C$.  Define $\Bound = \frac{1}{3} N \prod_{l
\mid N, l \text{ prime}} (1 + \frac{1}{l})$, so that the subgroup of
$\ZZ^{D \times D}$ (for $D$ the dimension of $C$) generated by
matrices representing the $T_n$ for $n \le \Bound$ equals the Hecke
algebra of weight $2$.  Let $\cA$ be the sub-$\FF_2$-vector space of
$\FF_2^{D\times D}$ generated by $\overline{T_n}$ for $n \le \Bound$.
Define $\cR$ to be the subspace of $\cA$ generated by $\overline{T_n}$
for all odd $n \le \Bound$.

Using the natural surjection $< T_i \; | \; i \le \Bound > \otimes_\ZZ
\FF_2 \twoheadrightarrow \cA,$ it follows immediately from the results
of the preceding section that the $\FF_2$-vector space
$$ \cH = \cA / \cR $$
is equipped with an action by the Hecke algebra of 
$\ksf 1 0 N {\FF_2}$ similar to the one explained in proposition~\ref{thm2}.

The function {\tt HeckeAlgebraWt1} of {\em Hecke1} computes this
module $\cH$ and also the first $\Bound$ Hecke operators of weight $1$
acting on it. More precisely, a record containing the necessary data
is created.  Properties can be accessed using e.g.~the commands {\tt
Dimension}, {\tt Field}, {\tt HeckeOperatorWt1}, {\tt HeckeAlgebra}
and {\tt HeckePropsToString}.  Please consult \cite{hdoc} (and
\cite{cma}) for a precise documentation of the provided functions.

\subsubsection*{An example session}

We assume that the packages {\em CommMatAlg} and {\em Hecke1} are
stored in the folder \ttt{PATH}. We attach the packages by typing
\comp{> Attach ("PATH/CommMatAlg.mg");\\
> Attach ("PATH/Hecke1.mg");}
We can now create a record containing all information for
computations of Hecke operators of weight $1$ acting on $\cH$ (as
described above with $N = 491$ and $p = 2$).
\comp{> M := ModularSymbols (491,2); \\
> h := HeckeAlgebraWt1 (M);}
It is not advisable to access information by printing {\tt h}.
Instead, we proceed as follows:
\comp{> Dimension(h);\\
6\\
> Bound(h);\\
164}
These functions have the obvious meanings. 
If one is interested in some properties of the Hecke algebra acting
on $\cH$, one can use:
\comp{> HeckePropsToString(h);\\[0mm]
Level N = 491:\\[0mm]
***************************\\[0mm]
\\[0mm]
Dimension = 6\\[0mm]
Bound = 164\\[0mm]
Class number of quadratic extension with |disc| = 491 is: 9\\[0mm]
There are 2 local factors.\\[0mm]
\\[0mm]
Looking at 1st local factor:\\[0mm]
\\[0mm]
  Residue field = GF(8)\\[0mm]
  Local dimension = 3\\[0mm]
  UPO = 1\\[0mm]
  Eigenvalues = \{ $1, w, w^2, w^4, 0$ \}\\[0mm]
  Number of max. ideals over residue field = 3\\[0mm]
\\[0mm]
Looking at 2nd local factor:\\[0mm]
\\[0mm]
  Residue field = GF(2)\\[0mm]
  Local dimension = 3\\[0mm]
  UPO = 3\\[0mm]
  Eigenvalues = \{ $0, 1$ \}\\[0mm]
  Number of max. ideals over residue field = 1}
Here $w$ stands for a generator of the residue field in question.
For the significance of these data, please see the following section.

\subsection{Mestre's calculations}

In this section we report on computations we performed in relation
with Mestre's calculations exposed in appendix A. Mestre considered
weight $1$ modular forms for $\Gamma_0(N)$, where $N$ is an odd prime.

According to the modified version of Serre's conjecture (see
e.g.~\cite{Edix5}), one expects that for any $2$-dimensional
irreducible Galois representation
$$ \rho: G_\QQ \to \SL_2(\overline{\FF_2}), $$
which is {\em unramified at $2$}, there exists a {\em weight $1$}
Hecke eigenform $f \in \ksf 1 0 {N_\rho} {\overline{\FF_2}}$ giving rise
to the representation $\rho$ via Deligne's theorem. Here $N_\rho$ is the
Artin conductor of the representation $\rho$.

Unfortunately, the implication $\rho$ is modular, hence $\rho$ comes
from a form of weight $1$ and level $N_\rho$ is unproved in the
exceptional case $p = 2$.

There is a simple way to produce Galois representations, which are
unramified at $2$, with given Artin conductor $N$, when $N$ is odd and
square-free.  One considers the quadratic field $K = \QQ(\sqrt{N})$
resp.  $K = \QQ(\sqrt{-N})$ if $N \equiv 1 \; (4)$ resp. $N \equiv 3
\; (4)$, which has discriminant $(N)$. Let now $L$ be the maximal
subfield of the Hilbert class field of $K$ such that $[L:K]$ is odd.
Then $L$ is Galois over $\QQ$ of degree $2u$ with $u$ the odd part of
the class number of $K$. The Galois groups in question form a split
exact sequence $ 0 \to G_{L|K} \to G_{L|\QQ} \to G_{K|\QQ} \to 0.$ The
conjugation action of $G_{K | \QQ}$ via the split on $G_{L|K}$ is by
inversion.  For any character $\chi: G_{L|K} \to \overline{\FF_2}^*$,
one has the induced representation $\Ind_{G_{L|K}}^{G_{L|\QQ}}(\chi):
G_{L|\QQ} \to \SL_2(\overline{\FF_2})$.  It is irreducible if $\chi$
is non-trivial, and $\Ind_{G_{L|K}}^{G_{L|\QQ}} (\chi_1) \cong
\Ind_{G_{L|K}}^{G_{L|\QQ}} (\chi_2)$ if and only if $\chi_1 = \chi_2$
or $\chi_1 = \chi_2^{-1}$.  The Artin conductor of
$\Ind_{G_{L|K}}^{G_{L|\QQ}}(\chi)$ is $N$.  Consequently, one receives
$(u-1)/2$ non-isomorphic Galois representations with dihedral image
and Artin conductor $N$.  More precisely, the image of
$\Ind_{G_{L|K}}^{G_{L|\QQ}}(\chi)$ is $D_{2\cdot
\#\Image(\chi)}$. These are the dihedral representations to which
Mestre refers in appendix~A.

It is known that any dihedral representation
$\Ind_{G_K}^{G_\QQ}(\chi)$ is modular, where $K|\QQ$ is a quadratic
field and $\chi: G_K \to \overline{\FF_2}^*$ is a character.  However,
as mentioned above, the weight and the level are not known to occur as
predicted.  Looking at the standard proof (see e.g. \cite{Fermat},
Theorem 3.14) of modularity, we see that obstacles occur if $K$ is
real and does not allow any non-real unramified quadratic extension.

A feature of modular forms over fields of positive characteristic is
that even for prime levels the Hecke algebra can be non-reduced.  The
Hecke algebra is finite-dimensional and commutative, hence it splits
into a direct product $\TT = \prod_{i=1}^r \TT_i$ of local algebras.
For a local algebra $\TT_i$ with maximal ideal $\fm_i$, we introduce
the number $\omega(\TT_i) = \min \{ \; n \; | \; (\fm_i)^n = (0) \;
\}$.  It is related to the number $B(m)$ considered in appendix A: one
is in the case $B(m)$ with $m \le \omega(\TT) := \max_i
(\omega(\TT_i))$.

Mestre considered all prime levels up to $1429$ and some higher ones.
The dimension we find for the space $(\TT(1) \otimes_\ZZ \FF_2)/
\tilde{\cR}$ (see Prop.~\ref{thm2}) equals the dimension announced by
Mestre.  Moreover, he finds case $B(m)$ if and only if we find $m =
\omega(\TT)$ (from the definition of the two numbers, the equality
does not follow in general).  We also calculated the image of the
Galois representations associated to the eigenforms we found. These
images agree with Mestre's claims.  More precisely, we compute that
for prime level $N$ less than $2100$ there exists an eigenform with
image equal to $A_5 = \SL_2(\FF_4)$ in the cases
$$ N = 653, 1061, 1381, 1553, 1733, 2029 $$
and equal to $\SL_2(\FF_8)$ in the cases
$$ N = 1429, 1567, 1613, 1693, 1997, 2017, 2089.$$ In all other prime
cases, we find only dihedral images. However, we always find all the
dihedral eigenforms predicted by the modified version of Serre's
conjecture.

One of the main points of Mestre's letter to Serre was to conclude
from the existence of an $\SL_2(\FF_8)$-form that not all weight $1$
forms arise as reductions of weight $1$ forms from characteristic $0$,
even for an increased level because $\SL_2(\FF_8)$ is not a quotient
of a finite subgroup of $\PGL_2(\CC)$.  We can reformulate that by
saying that whenever there is an $\SL_2(\FF_8)$-form, the space of
Katz modular forms of weight $1$ is strictly bigger than the space of
classical forms.

To finish with, we wish to point out that in prime levels the
representations associated to eigenforms of weight $1$ in
characteristic $2$ were always found to be irreducible and the Hecke
algebra to be of type Gorenstein.  For non-prime square-free levels
both properties can fail.


\begin{thebibliography}{99}
\bibitem{hdoc} G.\ Wiese. {\sl The Magma package Hecke1}, 
documentation and source are available on the author's homepage
\bibitem{cma} G.\ Wiese. {\sl The Magma package CommMatAlg}, 
documentation and source are available on the author's homepage
\end{thebibliography}
\end{document}